\documentclass[10pt,leqno]{article}
\usepackage{amsmath}
\usepackage{amsfonts}
\usepackage{mathrsfs}
\newtheorem{theorem}{\bf Theorem}[section]

\newtheorem{corollary}{\bf Corollary}[section]
\newtheorem{remark}{\bf Remark}[section]
\newtheorem{Lemma}{\bf Lemma}[section]

\newcommand{\R}{{I\!\!R}}

\newcommand{\bv}{{\bf v}}
\newcommand{\bz}{{\bf z}}
\newcommand{\be}{{\bf e}}
\newcommand{\bs}{{\boldsymbol \sigma}}
\newcommand{\bq}{{\bf q}}
\newcommand{\bp}{{\bf p}}
\newcommand{\bpsi}{{\boldsymbol\psi}}

\newcommand{\eq}{{\boldsymbol\eta_{\bf q}}}
\newcommand{\tq}{{\boldsymbol\theta_{\bf q}}}
\newcommand{\xq}{{\boldsymbol\xi_{\bf q}}}
\newcommand{\roq}{{\boldsymbol\rho_{\bf q}}}
\newcommand{\eu}{{\eta_u}}
\newcommand{\ut}{{\theta_u}}
\newcommand{\xu}{{\xi_u}}
\newcommand{\rou}{{\rho_u}}
\newcommand{\es}{{\boldsymbol\eta_{\boldsymbol\sigma}}}
\newcommand{\st}{{\boldsymbol\theta_{\boldsymbol\sigma}}}

\newcommand{\xs}{{\boldsymbol\xi_{\boldsymbol\sigma}}}

\newcommand{\eqt}{{\boldsymbol\eta_{\bq_t}}}

\newcommand{\xqt}{{\boldsymbol\xi_{{\bf q}_t} }}

\newcommand{\eut}{{\eta_{u_t}}}
 
\newcommand{\xut}{{\xi_{u_t} }}

\newcommand{\est}{{\boldsymbol\eta_{{\boldsymbol\sigma}_t}}}

\newcommand{\bxq}{\bar{\boldsymbol\xi}_{\bf q}}
\newcommand{\bxs}{\bar{\boldsymbol\xi}_{\boldsymbol\sigma}}

\newcommand{\bV}{{\bf V}}
\newcommand{\bH}{{\bf H}}

\newcommand{\iph}{{i+1/2}}

\newcommand{\xQh}{\hat{\boldsymbol\xi}_\Q}
\newcommand{\xZh}{\hat{\boldsymbol\xi}_\Z}

\newcommand{\eU}{\eta_U}
\newcommand{\eQ}{\boldsymbol\eta_\Q}
\newcommand{\eZ}{\boldsymbol\eta_\Z}
\newcommand{\xU}{\xi_U}
\newcommand{\xQ}{\boldsymbol\xi_\Q}
\newcommand{\xZ}{\boldsymbol\xi_\Z}

\newcommand{\ta}{\tau^\ast}

\newcommand{\B}{\mathcal {B}}

\newcommand{\ii}{{\boldsymbol \epsilon}}
\newcommand{\Bnj}{{\tilde B}^\nph_\jph}

\newcommand{\ph}{{\partial}_t}
\newcommand{\pc}{{\delta}_t}

\newcommand{\pb}{\bar{\partial}_t}
\newcommand{\ps}{{\partial}_{t}^2}
\newcommand{\half}{{1/2}}

\newcommand{\Z}{{\boldsymbol Z}}
\newcommand{\Q}{{\boldsymbol Q}}
\newcommand{\U}{{U}}
\newcommand{\dt}{k}
\newcommand{\qr}{1/4}

\newcommand{\nph}{{n+1/2}}
\newcommand{\msph}{{m^{\star}+1/2}}
\newcommand{\jph}{{j+1/2}}
\newcommand{\nmh}{{n-1/2}}

\newcommand{\ah}{A^\half}

\newcommand{\ce}{{\mathcal E}_{\B}}
\newcommand{\hce}{{\mathcal {\hat E}}_{\B}}
\newcommand{\Hd}{\bH(div;\Omega)}

\newcommand{\se}{\setcounter{equation}{0}}
\newsavebox{\savepar}

\begin{document}

\title{\bf Optimal error estimates of  mixed FEMs for second order
hyperbolic integro-differential equations with minimal smoothness on initial data}

\author{Samir Karaa{\footnote{
Department of Mathematics and Statistics, Sultan Qaboos University, 
P. O. Box 36, Al-Khod 123,
Muscat, Oman. Email: skaraa@squ.edu.om} } 
 and
Amiya K. Pani{\footnote{
Department of Mathematics, Industrial Mathematics Group, Indian
Institute of Technology Bombay,
Powai, Mumbai-400076. Email: akp@math.iitb.ac.in}}}

\date{}
\maketitle

\begin{abstract}
In this article, mixed finite element methods are discussed for a class 
of hyperbolic integro-differential equations (HIDEs). Based on a modification of the nonstandard energy formulation of Baker, both semidiscrete and completely discrete implicit schemes for an extended mixed method are analyzed and optimal  $L^\infty(L^2)$-error estimates are derived under minimal smoothness assumptions on the
initial data. Further, quasi-optimal estimates are shown to hold in $L^\infty(L^\infty)$-norm. 
Finally, the analysis is  extended to the standard mixed method for HIDEs and   optimal error  
estimates in $L^{\infty}(L^2)$-norm are derived again under minimal smoothness on initial data.
\end{abstract}


\section{\bf Introduction} 
We consider two mixed finite element methods for the following
hyperbolic  integro-differential equation:
\begin {eqnarray}
 u_{tt}-\nabla\cdot{\Big (}A\nabla u -\int_{0}^{t} B(t,s)\nabla u(s)ds{\Big )} & = & 0
 ~~~~~~~~~~\mbox{in}~~\Omega\times J,\label{1}\\
 u(x,t) & = & 0~~~~~~~~~~\mbox{on}~\partial \Omega\times J,\label{2}\\
 u(x,0) & = & u_0 ~~~~~\mbox{in}~\Omega, \label{3}\\
 u_t(x,0) & = & u_1 ~~~~~\mbox{in}~\Omega, \label{4}
\end{eqnarray}
with given functions $u_0$ and $u_1$, where $\Omega\subset \R^2$ is a bounded
convex polygonal domain, $J=(0,T],~T<\infty$, and $u_{tt}=\frac{\partial^2 u}
{\partial t^2}.$ Here, $A=[{a_{ij}(x)}]_{1\le i,j\le d}$ and 
$B(t,s)=[{b_{ij}(x,t,s)}]_{1\le i,j\le d}$ are $2\times 2$ matrices 
with smooth coefficients. Further, we
assume that $A$ is symmetric and uniformly positive definite in $\bar{\Omega}$. 
Problems of this kind arise in linear viscoelasticity specially in viscoelastic materials with memory 
(cf. Renardy  {\it et al.} \cite{RHN}).

Early {\it a priori} error estimates for Galerkin finite element methods for solving 
(\ref{1})-(\ref{4}) (without integral terms) 
were derived by Dupont \cite{Dupont-15} using a standard energy argument. 
These estimates were improved by Baker \cite{Baker-3}, who used a technique 
interpreted later in literature as a nonstandard energy arguments. In \cite{Rauch}, Rauch discussed 
the convergence of a continuous time Galerkin approximation to a second order
wave equation and proved optimal error estimates
in $L^\infty(L^2)$-norm using piecewise linear polynomials, when $u_0\in H^3\cap H_0^1$ (with $u_1=0$). The
analysis improves upon the earlier results of Baker and Dougalis \cite{BD-1996},
 where optimal error estimates 
were derived under the assumptions that $u_0\in H^5\cap H_0^1$ and $u_1\in H^4\cap H_0^1$.
In \cite{SP},  Sinha and  Pani extended Rauch's results to hyperbolic integro-differential
equation with quadrature and obtained  optimal error estimates
in $L^\infty(L^2)$ norm in the case $u_0\in H^3\cap H_0^1$ and $u_1\in H^2\cap H_0^1$. For more on
Galerkin methods and optimal error estimates for the problem (\ref{1})-(\ref{4}), see, 
\cite{CL-89,PTW}. 

In the  literature, optimal error estimates for mixed finite element approximations 
to second order hyperbolic equations were proved in 
\cite{CDW-13,Geveci-19,JRW,BJT,CE}. In \cite{Geveci-19}, Geveci
derived $L^\infty$-in-time, $L^2$-in-space error bounds for the
continuous-in-time mixed finite element approximations of velocity and stress.
In \cite{CDW-13}, {\it a priori} error estimates were obtained for mixed
finite element approximations to both displacement and stress requiring less regularity 
than needed in \cite{Geveci-19}. Stability for a family of discrete-in-time
schemes was also demonstrated under regularity $u_0\in H^4\cap H^1_0$ and $u_1 \in H^3,$ when $r=2.$ In \cite{JRW}, a mixed finite element displacement formulation was proposed 
for  the acoustic wave equation under reduced regularity  requirement on the displacement variable. 
Further, for a completely discrete  explicit scheme,
stability results and error estimates were established.

For the standard mixed methods, the problem (\ref{1})
is often rewritten by introducing a new variable
\begin{eqnarray}\label{sigma-0}
\bs(t)= A \nabla u - \int_{0}^{t} B(t,s) \nabla u(s)\,ds,
\end{eqnarray}
or equivalently
\begin{eqnarray}\label{sigma-1}
\alpha \bs(t)= \nabla u - \int_{0}^{t} A^{-1} B(t,s) \nabla u(s)\,ds,
\end{eqnarray}
as
\begin{eqnarray}\label{mixed-form}
u_{tt} -\nabla\cdot \bs(t)=0,
\end{eqnarray}
where $\alpha= A^{-1}.$ Further, using resolvent operator, 
(\ref{sigma-1}) is formulated as
\begin{eqnarray}\label{sigma-2}
\nabla u(t) = \alpha \bs (t) + \int_{0}^{t} M(t,s) \bs (s)\,ds,
\end{eqnarray}
where $M(t,s)= R(t,s) A^{-1} $ and $R(t,s)$ is the resolvent of the matrix
$A^{-1} B(t,s)$ given by 
\begin{eqnarray}\label{resolvent}
R(t,s)= A^{-1} B(t,s) +  \int_{s}^{t} A^{-1} B(t,\tau)R(\tau, s)
\,d\tau,~~~t> s \geq 0.
\end{eqnarray}
With $W=L^2(\Omega)$ and $\bV= \Hd,$ the weak formulation for
the mixed problem (\ref{mixed-form})-(\ref{sigma-2}) is to
seek a pair of  functions $(u, \bs):(0,T]\longrightarrow W \times \bV$ 
satisfying
\begin{eqnarray}\label{weak-mixed-1}
(\alpha \bs, \bv) + \int_{0}^{t} (M(t,s) \bs (s),\bv)\,ds + (\nabla\cdot \bv, u)
&=& 0, \;\;\;\bv \in \bV,\\ \label{weak-mixed-2}
(u_{tt}, w)- (\nabla\cdot \bs, w) &=& 0, \;\;\;w \in W
\end{eqnarray}
with $u(0)= u_0$  and $u_t(0)= u_1$. To the best of our knowledge, there are few 
 results available on optimal estimates of  mixed approximations to problem 
(\ref{1})-(\ref{4}) with minimum smoothness on initial data.  
In the context of  parabolic problems, based on  mixed methods related to the weak formulation
(\ref{weak-mixed-1})-(\ref{weak-mixed-2}), 
Sinha {\it et al.} \cite{SEL} have derived an optimal convergence rate $O(ht^{-1})$ 
 for the velocity $\bs(t)$ in ${\bf {L}}^2$-norm  and suboptimal convergence
rate $O(ht^{-1/2})$ for the pressure $u(t)$ in $L^2 $-norm with $t\in (0,T],$ 
when $u_0 \in L^2(\Omega).$ However, optimal rate 
$O(h^2 t^{-1})$ for $u$ is only established for a class of  problems when $A= aI$ 
and $B= b(t,s)I,$ where $a$ and $b$ are independent of spatial variable $x$. 
These results have recently been improved by Goswami {\it et al.}  \cite{GPY-2013}, 
who have established optimal convergence rate  $O(h^2 t^{-1})$ for $u(t)$ in  $L^2$-norm
and  $O(h t^{-1})$ for $\bs(t)$ ${\bf {L}}^2$-norm for all $t\in (0,T],$ when 
$u_0\in L^2(\Omega).$ Similar results have been obtained for the extended finite element 
method described below. In both papers only semidiscrete problems have been considered.

In the first part of this paper, an extended mixed method for (\ref{1}) is proposed and 
analyzed. Such a method was analysed earlier in \cite{Chen-1} and \cite{AWY} for elliptic problem, \cite{WD} for degenerate nonlinear parabolic problem, \cite{GPY-2013} for parabolic integro-differential equations and references cited in.
To motivate this new mixed method, we introduce  two variables:
\begin{eqnarray}\label{f1}
\bq=\nabla u,~~~~{\mbox{ and }}~~~~ \bs= A\bq-\int_0^t B(t,s)\bq(s)ds.
\end{eqnarray}
Then, the equation (\ref{1}) takes the form
 $$u_{tt}-\nabla\cdot\bs = 0.$$ 
Now, the weak mixed formulation of (\ref{1})-(\ref{4}) which forms
a basis of our mixed Galerkin method is to  
find $(u,\bq,\bs):J\rightarrow W\times\bV\times\bV$ satisfying
\begin{eqnarray}
(\bq,\bv)+(u,\nabla\cdot\bv)=0~~~\bv\in\bV,\label{w1}\\
(\bs,\bz)-(A\bq,\bz)+\int_0^t(B(t,s)\bq(s),\bz)ds=0~~~\bz\in\bV,\label{w2}\\
 (u_{tt},w)-(\nabla\cdot\bs,w) = 0~~~w\in W, \label{w3}
\end{eqnarray}
with $u(0)=u_0$ and $u_t(0)=u_1$.  The main goal of this paper is to 
establish optimal convergence rate for the approximate solutions of (\ref{1}),
when the initial functions $u_0\in H^3\cap H_0^1$ and $u_1\in H^2\cap H_0^1$.
Our analysis is essentially based  on a simple energy technique and a use of a time integration without exploiting the inverse of the associated discrete elliptic operator. Essentially, integrating in time
leads to a first order evolution process and hence, is instrumental in reducing the regularity 
requirements on the solution. Further, due to the presence of the integral term, it is observed  that the concept of 
mixed Ritz-Volterra projections used earlier  in \cite{ELSWZ},\cite{ELWZ05},\cite{SEL} and \cite{GPY-2013}
plays a crucial role in our analysis. For the completely discrete scheme, a major
difficulty associated with a use of the nonstandard energy formulation by Baker is the presence 
of a fourth order time derivative of the displacement $u$ in the bounds, see \cite{CDW-13}. Moreover, 
additional difficulty is caused by the quadrature error associated with 
a second order midpoint quadrature rule which is used to approximate the integral term. Special care is
needed to arrest these issues otherwise  for optimal convergence rate, we land up with 
higher regularity assumption on the initial data, namely; $u_0\in H^4\cap H_0^1$. Therefore, 
a modification of Baker's approach is  adopted and optimal error estimates in $\ell^{\infty}(L^2)$-norm are  derived, when $u_0\in H^3\cap H_0^1$ and 
$u_1\in H^2\cap H_0^1$. This technique is proved to be powerful and can successfully 
be applied even to the problem in \cite{SP}. Finally, the analysis has been extended to the standard mixed method  corresponding to the formulation (\ref{weak-mixed-1})-(\ref{weak-mixed-2}) and error analysis has been
briefly discussed.

Now compared  to the (\ref{weak-mixed-1})-(\ref{weak-mixed-2}), the extended or expanded method corresponding to (\ref{w1})-(\ref{w3}) may have introduced one more variable leading to a computation of one more extra variable. However, in Section~2, it  is shown that it is possible to eliminate $\bq_h,$ which is an approximation of the gradient vector $\bq.$ Hence, both these schemes have almost comparable computational cost. Moreover, in the new formulation, we need not invert the coefficient matrix $A.$

Throughout this article, we denote by $C,$ a generic positive constant which may vary from 
context to context, and whenever there arises no confusion, we would denote $u(t)$ simply as $u$ 
for the sake of convenience.

An outline of the paper is as follows. In Section 2, we give some {\it a priori}
bounds and regularity results for (\ref{w1})-(\ref{w3}), and briefly present the finite
element approximation of the extended mixed formulation (\ref{w1})-(\ref{w3}).
In Section 3, the extended mixed Ritz-Volterra
projection is introduced and analyzed.  In Section 4,  error estimates for  Galerkin approximations of $u, \bs$ and $\bq$ 
for the semidiscrete problem are derived. The completely discrete problem is discussed in Section 5 and optimal error estimates in $\ell^{\infty}(L^2)$-norm are established.
Finally, in Section 6, an extention of our analysis to the mixed formulation
(\ref{weak-mixed-1})-(\ref{weak-mixed-2}) with minimal smoothness assumptions on the initial data 
is briefly discussed.


\section {\bf Extended Mixed Finite Element Method}
\se


For our analysis, we shall use the standard  notations for $L^2(\Omega),~H^1_0(\Omega),~ 
H^m(\Omega)$ and $\Hd$ spaces with their norms and seminorms.
To be more specific, $L^2(\Omega)$ is equipped with inner product 
$(\cdot ,\cdot)$ and norm $\|\cdot\|.$ The norm on $\Hd$ is given by
$$ \|\bv\|_{\bH(div;\Omega)}= (\|\bv\|^2+\|\nabla\cdot\bv\|^2)^{1/2}.$$
Further, the standard Sobolev space $H^m(\Omega)$ of order $m$ is equipped
with its norm $\|\cdot \|_{H^m(\Omega)},$ which we denote simply by $\|\cdot \|_{m}$.
Since the matrix $A$ in (\ref{1}) is positive definite, there exist 
positive constants $a_0$ and $a_1$ such that
\begin{equation}\label{apd}
 a_0 \|\bs\|\leq \|\bs\|_{A}\le a_1\|\bs\|,~~~\mbox{where}~~
 \|\bs\|^2_{A}:=(A\bs,\bs).
\end{equation}
Further, assume that all coefficients of $B$ and their derivatives are bounded  
in  their respective domain of definitions by the positive constant $a_1$.
Under our assumptions on the domain and  on the coefficient matrix 
$A,$ we note that the following elliptic regularity result holds: there exists a 
positive constant $C$  such that for  $\phi\in H^2\cap H_0^1$
\begin{equation}\label{er}
 \|\phi\|_2 \leq C\|\nabla\cdot(A\nabla\phi)\|.
\end{equation}
For our subsequent use, we state without proof {\it a priori} estimates for $u$, $\bq$ and $\bs$ 
satisfying (\ref{w1})-(\ref{w3}) under appropriate regularity conditions
on the initial data $u_0$ and $u_1$. For more details, we refer to  \cite{SP} and \cite{Sinha-2002}. 
\begin{Lemma}\label{lm1}
 Let $(u,\bq,\bs)$ satisfy $(\ref{w1})$-$(\ref{w3})$. Then, there is a positive constant $C$ such that 
the following regularity results hold:
$$
\|D_t^ju(t)||+||D_t^{j-1}u(t)\|_1+
\|D_t^{j-1}\bs(t)||+||D_t^{j-1}\bq(t)\| \leq
C(T)(\|u_0\|_j+\|u_1\|_{j-1}),\quad j=1,\cdots,4,
$$
and 
$$
\|D_t^ju(t)||_2\leq C(T)(\|u_0\|_{j+2}+\|u_1\|_{j+1}),\quad j=0,1,2,
$$
where $D_t^j=(\partial^j/\partial_t^j)$.
\end{Lemma}

Based on the mixed formulation (\ref{w1})-(\ref{w3}) for the problem (\ref{1})-(\ref{3}), 
we now introduce the extended mixed finite element Galerkin method.
Let $\mathcal{T}_h$ be a regular triangulation of $\Omega$ by triangles of diameter at most $h$. 
Let $ \bV_h \times W_h $ denote a pair of finite element spaces satisfying the following
conditions:
\begin{eqnarray*}
 &(i)& ~\nabla\cdot\bV_h \subset W_h,~\mbox{and} \\
 &(ii)& \mbox{ there exists a linear operator } \Pi_h :\bV\rightarrow \bV_h
\mbox{ such that } \nabla\cdot\Pi_h=P_h(\nabla\cdot),
\end{eqnarray*}
\hspace{2.2cm} where $P_h :W\rightarrow W_h$ is the $L^2$-projection defined by
$$(\phi-P_h\phi,~w_h)=0,~~~\forall~w_h\in W_h,~\phi\in W.$$
Further, assume that the finite element spaces satisfy the following
approximation properties:
\begin{equation}\label{f2}
\|\bs-\Pi_h\bs\| \le Ch^{r} \|\nabla\cdot\bs\|_{r-1},~~~~\|u-P_h u\| 
\le Ch^r\|u\|_r,~~r=1,2.
\end{equation} 
and on a quasi-uniform mesh,
\begin{equation}\label{2.3a}
\|u-P_h u\|_{L^\infty(\Omega)}\le Ch^r|\log h|^\half\|u\|_{r+1},~~r=1,2.
\end{equation} 
%
%
Although, we can have several choices for $\bV_h$ and $W_h$, here we  consider only
the Raviart-Thomas elements of order one \cite{BF}. 
Note that $P_h$ and $\Pi_h$ satisfy
\begin{equation}
 (\nabla\cdot(\bs-\Pi_h\bs),w_h)=0,~~w_h\in W_h;~~~(u-P_hu,\nabla\cdot\bv_h)=0,
 ~~\bv_h \in \bV_h.
\end{equation}
Now, the corresponding semidiscrete mixed finite element formulation is to seek a triplet 
$(u_h,\bq_h,\bs_h):(0,T]\longrightarrow W_h\times\bV_h\times \bV_h$ satisfying
\begin{eqnarray}
 (\bq_h,\bv_h)+(u_h,\nabla\cdot\bv_h)=0~~\forall~~\bv_h\in \bV_h,\label{sw1}\\
 (\bs_h,\bz_h)-(A\bq_h,\bz_h)+\int_0^t(B(t,s)\bq_h(s),\bz_h)ds=0~~\forall~~\bz_h
 \in\bV_h,\label{sw2}\\
 (u_{htt},w_h)-(\nabla\cdot\bs_h,w_h) = 0~~\forall ~~w_h\in W_h, \label{sw3}
\end{eqnarray}
with initial data $u_h(0)$ and $u_{ht}(0)$ to be defined later. 
Since $W_h$ and $\bV_h$ are finite dimensional spaces,
the discrete problem (\ref{sw1})-(\ref{sw3}) leads to a linear system 
consisting of differential, integral and algebraic equations. 
Let $\{v_i\}_{i=1}^{N_1}$  and $\{\bpsi\}_{i=1}^{N_2}$ be the basis functions of the finite
element spaces $W_h$ and $\bV_h$, respectively. Let
$$u_h(t)=\sum_{i=1}^{N_1}\alpha_i(t)v_i(x),\quad  
\bq_h(t)=\sum_{i=1}^{N_2}\beta_i(t){\bf \psi}_i(x),\quad 
\bs_h(t)=\sum_{i=1}^{N_2}\gamma_i(t){\bf \psi}_i(x),$$ 
and ${\bf\alpha}(t)=(\alpha_1(t),\alpha_2(t),\cdots,\alpha_{N_1}(t))^T$,
${\bf\beta}(t)=(\beta_1(t),\beta_2(t),\cdots,\beta_{N_2}(t))^T$, and 
${\bf\gamma}(t)=(\gamma_1(t),\gamma_2(t),\cdots,\gamma_{N_2}(t))^T$. 
By choosing the test functions $\bv_h={\bf \psi}_j(x)$ and $\bz_h={\bf \psi}_j(x)$, 
for $j=1,2,\cdots,N_2$, in (\ref{sw1}) and (\ref{sw2}), respectively, and $w_h=v_j(x)$ 
in (\ref{sw3}), we obtain the following system:
\begin{eqnarray}
{\mathbb D}{\bf \beta}(t)+{\mathbb C}{\bf\alpha}(t)=0,\label{sw1-v}\\
{\mathbb D}{\bf\gamma}(t)-{\mathbb A}{\bf\beta}(t)+\int_0^t{\mathbb B}(t,s){\bf\beta}(s)\,ds=0,\label{sw2-v}\\
{\mathbb D}\alpha''(t)-{\mathbb C}^T{\bf\gamma}(t)=0,\label{sw3-v}
\end{eqnarray} 
 with ${\bf\alpha}(0)$, ${\bf\alpha}'(0)$, ${\bf\beta}(0)$, and 
${\bf\gamma}(0)$ are given from  the initial data of the system 
(\ref{sw1})-(\ref{sw3}). The matrices in (\ref{sw1-v})-(\ref{sw3-v}) are defined
as follows
$${\mathbb D}=[({\bf \psi}_i,{\bf \psi}_j)]_{N_2\times N_2}, \qquad
{\mathbb A}=[(A{\bf \psi}_i,{\bf \psi}_j)]_{N_2\times N_2},$$ 
$${\mathbb B}(t,s)=[(B(t,s){\bf \psi}_i,{\bf \psi}_j)]_{N_2\times N_2},
\qquad{\mathbb C}=[(v_i,\mbox{div}{\bf \psi}_j)]_{N_1\times N_2}.$$
From (\ref{sw1-v}), we obtain ${\bf \beta}(t)=-{\mathbb D}^{-1}{\mathbb C}{\bf\alpha}(t)$. 
After  elimination of ${\bf \beta}(t)$, this can be seen a system of integro-differential
equation
\begin{eqnarray}
{\mathbb D}\alpha''(t)-{\mathbb C}^T{\bf\gamma}(t)&=&0,\label{sw3-vv}\\
{\mathbb D}{\bf\gamma}(t)+{\mathbb A}{\mathbb D}^{-1}{\mathbb C}{\bf\alpha}(t)&=&
\int_0^t{\mathbb B}(t,s){\mathbb D}^{-1}{\mathbb C}{\bf\alpha}(s)\,ds=0.\label{sw2-vv}
\end{eqnarray} 
Using Picard's method, it is easy to check that the system (\ref{sw2-vv})-(\ref{sw3-vv})
has a unique solution. 

Compared to the two field formulation which is based on the mixed weak form (\ref{weak-mixed-1})-(\ref{weak-mixed-2}) and is stated in Section 6, the new mixed system (\ref{sw1})-(\ref{sw3}) which depends on three field formulation has one more additional 
 variable to compute. However, one variable, say ${\bf {q}}_h$ can be easily eliminated with negligible computational cost  and hence,  it can have comparable computational cost.  Further, there is 
 no need of inverting the coefficient matrix $A.$



\section {\bf Mixed Ritz-Volterra Type Projections}
\se

In this section, the extended mixed Ritz-Volterra projections are introduced and analyzed. 
The projections are defined  as follows:
Given $(u(t),~\bq(t),~\bs(t)) \in W \times \bV \times \bV,$ for 
$ t \in (0,T],$ find $(\tilde u_h ,\tilde\bq_h,\tilde\bs_h): (0,T]
\longrightarrow W_h \times \bV_h\times \bV_h $ satisfying 
\begin{eqnarray}
 (\eq,\bv_h)+(\eu,\nabla\cdot\bv_h)=0,~~~\bv_h\in\bV_h, \label{eu1} \\
 (\es,\bz_h)-(A\eq,\bz_h)+\int_0^t(B(t,s)\eq(s),\bz_h)ds=0,~~~\bz_h\in\bV_h,
 \label{eu2} \\
 (\nabla\cdot\es,w_h)=0,~~~w_h\in W_h, \label{eu3}
\end{eqnarray} 
where $\eu= (u-\tilde u_h),\;\eq= (\bq-\tilde\bq_h)$ and 
$\es= (\bs-\tilde\bs_h).$  Since $W_h$ and $\bV_h$ are finite dimensional 
spaces, the discrete problem (\ref{eu1})-(\ref{eu3}), for a given triplet $\{u, \bq,\bs\},$ 
leads to a system of linear equations combined with algebraic constraints
for $\{\tilde u_h,\tilde\bq_h, \tilde\bs_h\}.$ Note that  when $B=0,$ the system
has a unique solution, see \cite{Chen-1}. Now using theory of linear 
Volterra equations of second kind and Picard's iteration, it is 
straightforward to prove that the system (\ref{eu1})-(\ref{eu3}), for a 
given triplet $\{u, \bq,\bs\},$ has a unique solution  $\{\tilde u_h,
\tilde\bq_h, \tilde\bs_h\}.$   

In this section, we  discuss estimates of $\eu,\eq $ and $\es.$ 
Using definitions of $P_h$ and $\Pi_h,$ we rewrite $\eu,\eq$ and $\es$ as
\begin{eqnarray*}
 \eq&=&(\bq-P_h\bq)-(\tilde\bq_h-P_h\bq)=:\tq-\roq, \\
 \eu&=&(u-P_h u)-(\tilde u_h-P_h u)=:\ut-\rou, \\
 \es&=&\bs-\tilde\bs_h=\bs-\Pi_h\bs=:\st.
\end{eqnarray*}
Because the estimates of $\ut,\tq$ and $\st$ are known, it is sufficient to 
estimate $\rou,\roq.$ Now rewrite (\ref{eu1})-(\ref{eu3}) as
\begin{eqnarray}
 (\roq,\bv_h)+(\rou,\nabla\cdot\bv_h)&=&0 ~~~\forall~\bv_h\in\bV_h, \label{r1} \\
 ~~~ -(A\roq,\bz_h)+\int_0^t(B(t,s)\roq(s),\bz_h)\,ds &=&(\st,\bz_h)-(A\tq,\bz_h)
 \label{r2} \\
 &+& \int_0^t(B(t,s)\tq(s),\bz_h)\,ds ~~~\forall~\bz_h\in\bV_h, \nonumber \\
 (\nabla\cdot\st,w_h)& =& 0 ~~~\forall~w_h\in W_h. \label{r3}
\end{eqnarray}
Below, we derive estimates of $\es$ and $\eq$.  The analysis of  Section 4 of \cite{GPY-2013}
can be suitably modified to prove the following results in Lemmas~\ref{lm2} and \ref{lm3}, but for completeness,
we indicate the proofs here.
\begin{Lemma}\label{lm2} 
Let $(\eu,\eq,\es)$ be such that the system $(\ref{eu1})$-$(\ref{eu3})$
is satisfied. Then, there exists a constant C independent of $h$ such that
for $t\in (0,T]$ 
\begin{equation}\label{nq}
\|D_t^j\es(t) \|+\|D_t^j\eq(t) \|\leq C h^r \left(\|u_0\|_{j+r+1}+
\|u_1\|_{j+r}\right),\quad j=0,1,2,  \quad r=1,2,
\end{equation}
and
\begin{equation}\label{hnq}
  \|\eq(t)\|_{(\Hd)^*} \leq Ch^2(\|u_0\|_2+\|u_1\|_1) +\|\eu\|,
\end{equation}
where $(\Hd)^*$ is the dual of $\Hd$.
\end{Lemma}
Proof. First, observe that
$$
\|D_t^j\es\|\leq Ch^r\|\nabla\cdot D_t^j\bs\|_{r-1}\leq Ch^r\|D_t^j u\|_{r+1},\quad j=0,1,2,\quad r=1,2,
$$
and
$$
\|\nabla\cdot D_t^j\es\|\leq Ch^r\|\nabla\cdot D_t^j\bs\|_r\leq Ch^r\| D_t^ju\|_{r+2}.
$$
Next, choose $\bz_h=\roq$ in (\ref{r2}) to obtain
\begin{eqnarray*}
 \|A^{1/2}\roq\|^2&=&-(\st,\roq)+(A\tq,\roq)+\int_0^t
 (B(t,s)\tq(s),\roq)ds \\
 && -\int_0^t(B(t,s)\roq(s),\roq(t))ds.
\end{eqnarray*}
Then, a use of the Cauchy-Schwarz inequality with the boundedness of $B$ 
and the positive definiteness property of $A$ yields
\begin{equation}\label{lm03}
 \|\roq\| \le C(T,a_1) \left(\|\st\|+\|\tq\|
+\int_0^t (\|\tq(s)\|+\|\roq(s)\|)\,ds\right).
\end{equation}
\noindent
Notice that by (\ref{f1}), it follows that 
\begin{eqnarray}\label{hqt}
 \|\tq\|~\le~Ch^r\|\bq\|_r~\le~Ch^r\|u\|_{r+1},
\;\;\;{\mbox {and}}\;\;  \|\st\|\leq Ch^r\|\nabla\cdot\bs\|_{r-1} \le Ch^r\|u\|_{r+1}.
\end{eqnarray}
A substitution of (\ref{hqt}) in (\ref{lm03}) with  Lemma~\ref{lm1} shows that
$$
\|\roq\| \le Ch^r(\|u_0\|_{r+1}+\|u_1\|_r)+C\int_0^t \|\roq(s)\|\,ds.
$$
An application of Gronwall's Lemma yields
$$
 \|\roq\| \le Ch^r(\|u_0\|_{r+1}+\|u_1\|_r). 
$$
A use of the triangle inequality establishes the estimate (\ref{nq}) for $j=0$.
Now,  differentiate (\ref{r2}) with respect to time to obtain
\begin{eqnarray*}
 (A\roq_{t},\bz_h)&=& -(\st_{t},\bz_h)+(A\tq_{t},\bz_h)-(B(t,t)\tq(t)
 ,\bz_h) \\
 &&-\int_0^t(B_t(t,s)\tq(s),\bz_h)\,ds+(B(t,t)\roq(t),\bz_h)+
 \int_0^t(B_t(t,s)\roq(s),\bz_h)~ds. 
\end{eqnarray*}
Again, apply the Cauchy-Schwarz inequality and the boundedness of $A$ and $B$ to arrive at
$$
 \|\roq_t\| \le C(T,a_1) \left(\|\roq\|+\|\tq\|+\|\tq_t\|+\|\st_t\|
+\int_0^t (\|\tq(s)\|+\|\roq(s)\|)\,ds\right).
$$
Taking into account approximation properties (\ref{f2}) and (\ref{f1}), it follows that
$$
 \|\tq_t\|~\le~Ch^r\|\bq_t\|_r~\le~Ch^r\|u_t\|_{r+1},
\;\;\;{\mbox {and}}\;\;  \|\st_t\| \le Ch^r\|u_t\|_{r+1},
$$
and hence, by Lemma~\ref{lm1}, 
$$
 \|\roq_t\| \le C(T)h^r(\|u_0\|_{r+2}+\|u_1\|_{r+1}). 
$$
Now, a use of the triangle inequality completes the proof of (\ref{nq}) for $j=1$. For $j=2$,
 differentiate again (\ref{r2}) with respect to $t$ and repeat the above arguments to obtain the estimate.

For the second estimate (\ref{hnq}), we use (\ref{eu1})  for any $\bv\in\bV$ to 
arrive at
\begin{eqnarray*} \label{eq-hat-hdiv}
 (\eq,\bv)&=&(\eq,\bv-\Pi_h\bv)+(\eq,\Pi_h\bv) \\
 &=&(\eq,\bv-\Pi_h\bv)-(\eu,\nabla\cdot \Pi_h\bv) \nonumber\\
 &\leq& \|\eq\|\|\bv-\Pi_h\bv\|+\|\eu\|\|\nabla\cdot \Pi_h\bv\| \nonumber\\
 &\le& \Big(Ch^2(\|u_0\|_2+\|u_1\|_1)+\|\eu\|\Big)\;\|\bv\|_{\bH(div;\Omega)},\nonumber
\end{eqnarray*}
and hence, for nonzero $\bv\in\bV$
$$
\frac{(\eq,\bv)}{\|\bv\|_{\bH(div;\Omega)}}\leq Ch^2(\|u_0\|_2+\|u_1\|_1)+\|\eu\|.
$$
By taking supremum over all nonzero $\bv\in\bV$, we obtain the desired
estimate and this  concludes the  proof. \hfill{$\Box$} \vspace{0.5cm}

Now, we  use a duality argument to estimate $\|\eu\|$.

\begin{Lemma}\label{lm3}
 Let $(\eu,\eq,\es)$ satisfy the system $(\ref{eu1})$-$(\ref{eu3})$. 
Then, there is a positive constant C independent of $h$ such that
for $t\in (0,T]$
 \begin{equation}\label{heu}
  \|D_t^j\eu(t)\| \leq Ch^2\left(\|u_0\|_{2+j}+\|u_1\|_{1+j}\right),\quad j=0,1,2,
 \end{equation}
and
\begin{equation}\label{3.15a}
\|\eu(t)\|_{L^\infty(\Omega)}\le
Ch^2|\log h|\left(\|u_0\|_{3}+\|u_1\|_{2}\right).
\end{equation}
\end{Lemma}
\noindent
Proof. Consider the auxiliary elliptic  problem:
\begin{eqnarray} 
 \nabla\cdot(A\nabla\zeta) &=& \eu ~\mbox{in}~\Omega, \label{aux01} \\
 \zeta &=& 0~~\mbox{on}~\partial\Omega, \label{aux02}
\end{eqnarray}
\noindent
and set $ \bp=\nabla\zeta$ and $\bpsi=A\bp$ so that $\nabla\cdot\bpsi=\eu.$
Then, from the elliptic regularity result (\ref{er}), it follows that
\begin{equation}\label{er1}
 \|\zeta\|_2,~\|\bp\|_1,~\|\bpsi\|_1 \le C\|\eu\|.
\end{equation}
Clearly, the following system of equations is satisfied for all $(w,\bv,\bz)\in W\times
\bV\times\bV$ 
\begin{eqnarray}
(\bp,\bv)+(\zeta,\nabla\cdot\bv)=0,\label{dh1}\\
(\bpsi,\bz)-(A\bp,\bz)=0,\label{dh2}\\
(\nabla\cdot\bpsi,w)=(\eu,w).\label{dh3}
\end{eqnarray}
Now, choose $w=\eu$, $\bz=\eq$ and $\bv=\es$ in (\ref{dh1})-(\ref{dh3}),
respectively,  and then add the resulting equations to arrive at
\begin{eqnarray}\label{heu-1}
 \|\eu\|^2 
 &=&(\nabla\cdot\bpsi,\eu)+(\bpsi,\eq)-(A\bp,\eq)+(\bp,\es)+(\zeta,
 \nabla\cdot\es) \\
 &=& (\nabla\cdot(\bpsi-\Pi_h\bpsi),\eu)+(\bpsi-\Pi_h\bpsi,\eq)-(A(\bp-\Pi_h\bp),\eq)\nonumber \\
 &&+ (\bp-\Pi_h\bp,\es)+(\zeta-P_h\zeta,\nabla\cdot\es)+(\nabla\cdot\Pi_h\bpsi,\eu) \nonumber\\
 &&+ (\Pi_h\bpsi,\eq)-( A(\Pi_h\bp),\eq)+(\Pi_h\bp,\es )+(P_h\zeta,\nabla\cdot\es).\nonumber
\end{eqnarray}
\noindent
Next, set $\bv_h=\Pi_h\bpsi$, $\bz_h=\Pi_h\bp$ and $w_h=P_h\zeta$ in 
(\ref{eu1})-(\ref{eu3}) with $\eu=\ut-\rou.$ 
Then substitute in (\ref{heu-1}) and use the Cauchy-Schwarz inequality and (\ref{f2}) 
for $r=1$ to obtain
\begin{eqnarray*}\label{heu-2}
 \|\eu\|^2&=&(\nabla\cdot(\bpsi-\Pi_h\bpsi),\ut)+(\bpsi-\Pi_h\bpsi,\eq)
 -(A(\bp-\Pi_h \bp),\eq)\\
 &&+(\bp-\Pi_h\bp,\es) +(\zeta-P_h\zeta,\nabla\cdot\es)-\int_0^t(B(t,s)
 \eq(s),\Pi_h\bp)\,ds \nonumber\\
 &\leq&
C\left(h^2\|u\|_2\|\nabla\cdot\bpsi\|+h\|\eq\|\|\nabla\cdot\bpsi\|+h\|\eq\|\|\nabla\cdot\bp\|\right)\nonumber\\
&&+C\left(h^2\|\nabla\cdot\bp\|\|\nabla\cdot\bs\|+h^2\|\zeta\|_2\|\nabla\cdot\bs\|\right)+
\int_0^t|(B(t,s)\eq(s),\Pi_h\bp)|\,ds. \nonumber
\end{eqnarray*}
Note that for fixed $s,t,$ we can use  (\ref{f2}) for $r=1$ to arrive at 
\begin{align}\label{B-term}
(B(t,s)\eq(s), & \Pi_h\bp)= -(B(t,s)\eq(s),\bp-\Pi_h\bp)  +(\eq(s),B^*(t,s) \bp) \nonumber\\
&\le C\|\eq(s)\|~\|\bp-\Pi_h \bp\|+\|\eq(s)\|_{(\bH(div;\Omega))^*} ~\|B^*(t,s) \bp\|_{\bH(div; \Omega)} \nonumber \\
&\le C\Big(h \|\eq(s)\| +\|\eq(s)\|_{(\bH(div;\Omega))^*}\Big)\,\|\bp\|_{1},   
\end{align}
Hence,
\begin{align}\label{hh2}
 \|\eu\|^2 
 &\le C \Big( h^2\|u\|_2 + h\|\eq\| + h^2\|\nabla\cdot\bs\| \Big)\, 
 \Big(\|\nabla\cdot\bpsi\|+\|\nabla\cdot\bp\| +  \|\zeta\|_2  \Big)  \nonumber\\
&+C\left(\int_0^t\left(h \|\eq(s)\| +
\|\eq(s)\|_{(\bH(div;\Omega))^*}\right)ds\right)\|\bp(t)\|_{1}. \nonumber
\end{align}
A use of estimates of $\|\eq\|$ from (\ref{nq}) with the elliptic regularity (\ref{er1}) yields
\begin{equation}\label{hh2}
\|\eu\|\le C h^2(\|u_0\|_2+\|u_1\|_1) + C\int_0^t \|\eq(s)\|_{(\bH(div;\Omega))^*}\,ds.
\end{equation}
Substituting (\ref{hnq}) in (\ref{hh2}), apply Gronwall's lemma to obtain
$$
 \|\eu\|\leq C h^2(\|u_0\|_2+\|u_1\|_1),
$$
and this concludes the proof of  (\ref{heu}) for $j=0$.

In order to estimate $\|\eut\|$,   consider again the elliptic
problem (\ref{aux01})-(\ref{aux02}) with replacing ${\eta}_u$ on 
the right hand side of (\ref{aux01}) by $\eut$. By setting
$$ \bp=\nabla \zeta,~~~\bpsi=A\bp, $$
we have  that
\begin{eqnarray*}\label{t1} 
(\bp,\bv)+(\zeta,\nabla\cdot\bv)=0,\;\; 
 (\bpsi,\bz)-(A\bp,\bz)=0, \;\; {\mbox { and }}\;
 (\nabla\cdot\bpsi,w)=(\eut,w). 
\end{eqnarray*}
From the standard regularity results, it follows that
\begin{equation}\label{er3}
 \|\zeta\|_2,~\|\bp\|_1,~\|\bpsi\|_1 \le C\|\eu_t\|.
\end{equation}
Now, differentiate with respect to time the three equations in (\ref{eu1})-(\ref{eu1})
 to obtain
\begin{eqnarray}
 &&~~~(\eqt,\bv_h)+(\eut,\nabla\cdot\bv_h)=0,~~~\bv_h\in\bV_h, 
\label{de1} \\
 &&~~~(\est,\bz_h)-(A\eqt,\bz_h)+(B(t,t)\eq,\bz_h)+\int_0^t (B_t(t,s)\eq(s),
 \bz_h)ds=0,~~~\bz_h\in\bV_h, \label{de2} \\
 &&~~~(\nabla\cdot\est,w_h)=0,~~~w_h\in W_h. \label{de3}
\end{eqnarray}
Following the previous steps for deriving the estimate of $\|\eu\|$, a use of (\ref{t1})-(\ref{de3})
at the appropriate steps leads to
\begin{eqnarray}
 \|\eut\|^2 
 &=& (\nabla\cdot(\bpsi-\Pi_h\bpsi),~\eut)+(\eqt,~\bpsi-\Pi_h\bpsi)
-(A(\bp-\Pi_h\bp), ~\eqt) \nonumber\\
&&+(\est,~\bp-\Pi_h\bp)-(\est,~\bp)-(B(t,t)\eq,~\Pi_h\bp)
-\int_0^t(B_t(t,s)\eq, ~\Pi_h\bp)\,ds \label{eta-ut-estimate} \\
 &=& (\nabla\cdot(\bpsi-\Pi_h\bpsi),~\theta_{ut})+(\eqt,~\bpsi-\Pi_h\bpsi)
-(A(\bp-\Pi_h\bp), ~\eqt) \nonumber \\
&&+(\est,~\bp-\Pi_h\bp)
+ (B(t,t)\eq,\bp-\Pi_h\bp)+\int_0^t (B_t(t,s)\eq,\bp-\Pi_h\bp)\,ds  \nonumber\\ 
&&+ (\zeta-P_h\zeta,~\nabla\cdot\est) -(B(t,t)\eq,\bp) - \int_0^t (B_t(t,s)\eq,\bp)\,ds, \nonumber
\end{eqnarray}
where $P_h$  is the $L^2$-projection onto the conforming finite element space consisting
of $C^0$-piecewise linear elements which is a subspace of $H^1_0(\Omega)$.
Since for all other terms except the last three terms on the right hand side 
of (\ref{eta-ut-estimate}), it is easy to derive the estimates, it is enough to
estimate the last three terms. An application of integration by parts shows that
$$ 
|(\zeta-P_h\zeta,~\nabla\cdot\est)|
=|(\nabla(\zeta-P_h\zeta),~\est)|\leq\|\nabla(\zeta-P_h\zeta)\|
   \|\est\|\leq Ch^2\|u_t\|_2 \|\zeta\|_2.
$$
For the last two terms on the right hand side of (\ref{eta-ut-estimate}),  we observe
that
$$
(B(t,t)\eq, \bp)\le M\|\eq\|_{(\bH(div;\Omega))^*}
\|\bp\|_{1},
$$
and a use of integration by parts yields
$$
\left|\int_0^t(B_t(t,s)\eq, \bp)ds\right|
\le M\left( \int_0^t \|\eq(s)\|_{(\bH(div;\Omega))^*}\,ds \right) 
\|\bp(t)\|_{1}. $$
All together, we obtain using the Cauchy-Schwarz inequality,
Lemmas \ref{lm1}, \ref{lm2} and \ref{lm3}, the approximation properties 
(\ref{f2}) of  projections $P_h$ and $\Pi_h,$
with elliptic regularity result (\ref{er3}) in (\ref{eta-ut-estimate})
the following estimate
$$
\|\eta_{u_t}\| \leq C h^2 (\|u_0\|_3+\|u_1\|_2).
$$
Finally, by differentiating (\ref{de1})-(\ref{de3}) with respect to time
and following the previous steps we establish the estimate for $j=2$. 

In order to show the estimate (\ref{3.15a}), we write (\ref{eu1}) as
$$(\eq,\bv_h)-(\rou,\nabla\cdot\bv_h)=0,~~~\bv_h\in\bV_h.$$
Using Lemma~2.1 in \cite{CT-81}, it follows that
$$ \|\rou\|_{L^\infty(\Omega)}\le C|\log h|\,\|\eq\|,$$
for some constant $C$ independent of $h$. A use of (\ref{2.3a}) and (\ref{nq})
 completes the proof.
\hfill{$\Box$}

\section {\bf Semidiscrete Error Estimates}
\se
In this section, error estimates for the  semidiscrete problem are derived.
Using the mixed Ritz-Volterra  projections defined in section 3, we rewrite 
\begin{eqnarray*}
e_u &:=& u-u_h=(u-\tilde u_h)-(u_h-\tilde u_h)=:\eu-\xu, \\
{\bf e}_\bq &:=& \bq-\bq_h=(\bq-\tilde\bq_h)-(\bq_h-\tilde\bq_h)=:\eq-\xq, \\
{\bf e}_\bs &:=& \bs-\bs_h=(\bs-\tilde\bs_h)-(\bs_h-\tilde\bs_h)=:\es-\xs,
\end{eqnarray*}
where, $(u,\bq,\bs)$ and $(u_h,\bq_h,\bs_h)$ are solutions  of 
(\ref{w1})-(\ref{w3}) and (\ref{sw1})-(\ref{sw3}), respectively. 
Note that, $(e_u,{\bf e}_\bq,{\bf e}_\bs)$ satisfy the following  equations
\begin{eqnarray}
 ({\bf e}_{\bq},\bv_h)+({\bf e}_u,\nabla\cdot\bv_h)=0,~~~\bv_h\in\bV_h, \label{e1} \\
 ({\bf e}_\bs,\bz_h)-(A{\bf e}_\bq,\bz_h)+\int_0^t(B(t,s){\bf e}_\bq(s),\bz_h)ds
 =0,~~~\bz_h\in\bV_h, \label{e2} \\
 ({\bf e}_{u_{tt}},w_h)-(\nabla\cdot{\bf e}_\bs,w_h) =0,~~~w_h\in W_h.\label{e3}
\end{eqnarray}
Since estimates of $\eu, \eq $ and $\es$ are known from Lemmas \ref{lm2} 
and \ref{lm3}, it is sufficient to estimate $\xu,\xq$ and $\xs.$ 
Using  (\ref{eu1})-(\ref{eu3}), we rewrite (\ref{e1})-(\ref{e3}) as
\begin{eqnarray}
 (\xq,\bv_h)+(\xu,\nabla\cdot\bv_h)=0,~~~\bv_h\in\bV_h, \label{xi1} \\
 (\xs,\bz_h)-(A\xq,\bz_h)+\int_0^t(B(t,s)\xq(s),\bz_h)ds=0,~~~\bz_h\in\bV_h,
 \label{xi2} \\
 (\xi_{u_{tt}},w_h)-(\nabla\cdot\xs,w_h) 
= (\eta_{u_{tt}},w_h),~~~w_h\in W_h. \label{xi3}
\end{eqnarray}
Below, one of the main results for the semidiscrete problem is proved.
\begin{theorem} \label{thm1} 
 Let $(u,\bq,\bs)$ and $(u_h,\bq_h,\bs_h)$ satisfy $(\ref{w1})$-$(\ref{w3})$ 
and $(\ref{sw1})$-$(\ref{sw3})$, respectively with $u_h(0)= P_h u_0$ and 
$u_{ht}(0)= P_h u_1$. Then,
 there exists a positive constant $C$ independent of the discretizing
parameter $h$ such that for $t\in (0,T]$
\begin{eqnarray}
 \|u_t(t)-u_{ht}(t)\| \le  Ch^2\left(\|u_0\|_4+\|u_1\|_3\right), \label{mt01} 
\end{eqnarray}
and 
\begin{eqnarray}
 \|\bq(t)-\bq_h(t)\| + \|\bs(t)-\bs_h(t)\| \le  C h^2 \left(\|u_0\|_4+\|u_1\|_3\right).
 \label{mt02}
\end{eqnarray}
\end{theorem}
Proof. First differentiate (\ref{xi1}) with respect to time and set $\bv_h=\xs$ 
in the resulting equation, $\bz_h=-\xqt$ in (\ref{xi2}) and $w_h=\xut$ in (\ref{xi3}). 
Then, add the resulting equations to arrive at
\begin{equation} \label{estimate-xit}
\frac{1}{2}\frac{d}{dt}\left(\|\xi_{u_t}\|^2 + ||\ah\xq (t)||^2\right)=
(\eta_{u_{tt}}, \xi_{u_t})+\int_0^t(B(t,s)\xq(s),\boldsymbol\xi_{\bq_t})\,ds.
\end{equation}
Apply integration by parts to the integral term on the right hand side of 
(\ref{estimate-xit}) to find that 
\begin{eqnarray} \label{integral-term}
\int_0^t B(t,s;\xq(s),\boldsymbol\xi_{\bq_t})\,ds &=&
\frac{d}{dt}\int_0^t(B(t,s)\xq(s),\xq)\,ds - (B(t,t)\xq(t),\xq) \\
&&-\int_0^t(B_t(t,s)\xq(s),\xq)\,ds. \nonumber
\end{eqnarray}
Substitute (\ref{integral-term}) in (\ref{estimate-xit}), then integrate 
the resulting equation from $0$ to $t$. 
A use of the Cauchy-Schwarz inequality with
the boundedness of $B$ and the positive definite property of $A$ yields 
\begin{eqnarray}
\|\xi_{u_t}(t)\|^2 &+&||\ah\xq (t)||^2 \leq \|\xi_{u_t}(0)\|^2+||\ah\xq (0)||^2
+2\int_0^t||\eta_{u_{tt}}||\,||\xi_{u_t}||\;ds\label{xi2b} \\
&&+C(T,a_0,a_1)\left(\int_0^t||\ah\xq(s)||\,||\ah\xq(t))||\;ds
+ \int_0^t||\ah\xq(s)||^2\,ds\right).\nonumber
\end{eqnarray}
For some $t^\ast\in[0,t]$, let
$$
\left(\|\xi_{u_t}(t^\ast)\|^2+||\ah\xq (t^\ast)||^2\right)=\max_{0\leq s\leq t}
\left(\|\xi_{u_t}(t)\|^2+||\ah\xq (t)||^2 \right).
$$
Then, at $t=t^\ast$, (\ref{xi2b}) becomes 
\begin{eqnarray*}
\|\xi_{u_t}(t^\ast)\|+||\ah\xq (t^\ast)||&\leq& \|\xi_{u_t}(0)\|+||\ah\xq (0)||
+2\int_0^{t^\ast}\|\eta_{u_{tt}}\|\;ds\\
&&+C(T,a_0,a_1)\int_0^{t^\ast}\|\ah\xq(s)\|\;ds,
\end{eqnarray*}
and hence,
\begin{eqnarray*}
\|\xi_{u_t}(t)\|&+&||\ah\xq (t)||\leq \|\xi_{u_t}(t^\ast)\|+||\ah\xq (t^\ast)||\\
&\leq &\|\xi_{u_t}(0)\|+||\ah\xq (0)||
+2\int_0^t\|\eta_{u_{tt}}\|\;ds
+C(T,a_0,a_1)\int_0^t \|\ah\xq(s)\|\;ds.
\end{eqnarray*}
Now an application of  Gronwall's lemma yields
\begin{equation} \label{estimate-xit2}
\|\xi_{u_t}(t)\|+||\ah\xq(t)||
\leq C \Big(\|\xi_{u_t}(0)\|+||\xq (0)||+\int_0^t||\eta_{u_{tt}}||ds\Big).
\end{equation}
To estimate $\|\xs\|$,  choose $\bz_h=\xs$ in (\ref{xi2}). Then, use the
Cauchy-Schwarz inequality to arrive at
$$
\|\xs(t)\|\leq C\Big(\|\ah\xq\|+\int_0^t\|\ah\xq(s)\|\,ds\Big).
$$
From the triangle inequality, we obtain
 $$
 \|e_{u_t}\|\le\|\xi_{u_t}\|+\|\eta_{u_t}\|.
$$
Apply  Lemmas~\ref{lm2} and \ref{lm3}  with the choices 
$u_h(0)=Pu_0$,  $u_{ht}(0)=P_hu_1$ and $\bq_h(0)=\Pi_h(\nabla u_0)$ to arrive at the
 estimate (\ref{mt01}). In a similar way, 
we can establish (\ref{mt02}) and this completes the rest of the proof. \hfill{$\Box$}

As a consequence of Theorem~$\ref{thm1}$, we have the  following $L^\infty(L^\infty)$
estimate. 
\begin{corollary} \label{co1} 
Assume that the mesh is  quasi-uniform. Then, under the assumptions 
of Theorem~$\ref{thm1},$  there exists a positive constant $C$ independent 
of the discretizing parameter $h$ such that for $t\in (0,T]$
\begin{equation}\label{eu2-1cc}
\|u(t)-u_{h}(t)\|_{L^\infty(\Omega)}\leq 
Ch^2|\log h|\left(\|u_0\|_{4}+\|u_1\|_{3}\right).
\end{equation}
\end{corollary}
Proof. Apply Lemma~2.1 in \cite{CT-81} to (\ref{xi1}) to obtain
$$\|\xu(t)\|_{L^\infty(\Omega)}\leq C|\log h|\,\|\xq\|.$$
Since
$$
\|u(t)-u_{h}(t)\|_{L^\infty(\Omega)}\leq \|\eu(t)\|_{L^\infty(\Omega)}
+\|\xu(t)\|_{L^\infty(\Omega)},
$$
(\ref{eu2-1cc}) follows from (\ref{3.15a}) and (\ref{estimate-xit2}) and this  
completes the proof. \hfill{$\Box$}


\begin{remark}
As a consequence of Theorem~$\ref{thm1}$ and the following inequality
\begin{eqnarray*}
\|\xi_{u}(t)\|\leq C\Big(\|\xi_{u}(0)\|+ \int_0^t\|\xi_{u_t}\|\,ds\Big),
\end{eqnarray*}
we easily derive the following estimate of $u-u_h$:
\begin{equation}\label{eu2-1}
\|u-u_{h}\|_{L^\infty(L^2(\Omega))}
\leq Ch^2\; \left(\|u_0\|_4+\|u_1\|_3\right).
\end{equation}
\end{remark}


Since from (\ref{eu2-1}), we obtain an optimal error estimate of $\|u-u_h\|_{L^{\infty}(L^2(\Omega))},$  when $u_0\in H^4\cap H_0^1$ and $u_1\in H^3\cap H_0^1$,
we now use a variant of Baker's nonstandard formulation (see \cite{Baker-3}) to
provide a proof of $L^{\infty}(L^2)$ estimate of $u-u_h$ under reduced regularity conditions 
on $u_0$ and $u_1$. More precisely, we shall obtain optimal $L^\infty(L^2)$ estimate for 
$u-u_h,$ when $u_0\in H^3\cap H_0^1$ and $u_1\in H^2\cap H_0^1$. In the rest of this section, we make
use of the following notation for $\bar{\phi}$:
$$
\bar{\phi}(t)=\int_0^t\phi(s)ds.
$$
After integrating (\ref{xi2}) and (\ref{xi3}) with respect to $t$, the following 
system of equations is derived
\begin{eqnarray}
 (\xq,\bv_h)+(\xu,\nabla\cdot\bv_h)=0,~~~\bv_h\in\bV_h, \label{xi1-n} \\
 (\bxs,\bz_h)-(A\bxq,\bz_h)+\int_0^t\Big(\int_0^s (B(s,\tau)\xq(\tau),\bz_h)d\tau\Big)\,ds=0,~~~\bz_h\in\bV_h,
 \label{xi2-n} \\
 (\xi_{u_{t}},w_h)-(\nabla\cdot\bxs,w_h) 
= (\eta_{u_{t}}, w_h),~~~w_h\in W_h. \label{xi3-n}
\end{eqnarray}
Note that in (\ref{xi3-n}), we have used the fact that $u_{ht}(0)= P_h u_1,$ that is,
$(e_{ht}(0), w_h)=0$ for all $w_h\in W_h.$
\begin{theorem} \label{thm2} 
 Let $(u,\bq,\bs)$ and $(u_h,\bq_h,\bs_h)$ satisfy $(\ref{w1})$-$(\ref{w3})$ 
and $(\ref{sw1})$-$(\ref{sw3})$, respectively, with $u_h(0)= P_h u_0$ and 
$u_{ht}(0)= P_h u_1$. Then,
 there exists a positive constant $C$ independent of the discretizing
parameter $h$ such that for $t\in (0,T]$
\begin{equation}\label{eu3-n}
 \|u(t)-u_h(t)\| \le  Ch^2\; \left(\|u_0\|_3+\|u_1\|_2\right). 
\end{equation} 
\end{theorem}
Proof.
Choose  $\bv_h=\bxs$, ${\bz}_h=-\xq$ and $w_h=\xi_{u}$, respectively in (\ref{xi1-n}), 
(\ref{xi2-n}) and (\ref{xi3-n}). On adding the resulting equations, we arrive at 
$$
\frac{1}{2}\frac{d}{dt}\left[\|\xi_{u}\|^2+\|\ah\bxq\|^2\right]=(\eta_{u_{t}}, \xu)+
\int_0^t\Big(\int_0^s(B(s,\tau)\xq(\tau),\xq)d\tau\Big)\,ds.
$$
Integrate from $0$ to $t$ to deduce 
\begin{eqnarray}\label{xiu-1}
\|\xi_{u}(t)\|^2&+& \|A^\half\bxq(t)\|^2= \|\xi_{u}(0)\|^2+2\int_0^t(\eta_{u_{t}},\xu)\,ds\\
&&-2\int_0^t\int_0^s \int_0^\tau (B(\tau,\ta)\xq(\ta),\xq(s))\;d\ta \;d\tau \,ds.\nonumber
\end{eqnarray}
Let $I$ denote the last term on the right hand side of (\ref{xiu-1}).   
Integration by parts yields
\begin{eqnarray*}
I &=& 
-2\int_0^t\int_0^s (B(\tau,\tau)\bxq(\tau),\xq(s))\;d\tau ds +
2\int_0^t\int_0^s \int_0^\tau (B_{\ta}(\tau,\ta)\bxq(\ta),\xq(s))\,d\ta \,d\tau ds\\
&=&-2I_1+2I_2.
\end{eqnarray*}
For $I_1$,  we again integrate by parts in time so that
\begin{eqnarray*}
|I_1|&=&\left|\int_0^t (B(s,s)\bxq(s),\bxq(t))ds-\int_0^t
(B(s,s)\bxq(s),\bxq(s))ds\right|\\
&\leq & \frac{a_1}{a_0}\left\{\|A^\half\bxq(t)\|\int_0^t\|A^\half\bxq(s)\|ds+
\int_0^t\|A^\half\bxq(s)\|^2ds\right\}.
\end{eqnarray*}
Similarly for $I_2$, we note that
\begin{eqnarray*}
|I_2|&=& \left|\int_0^t\int_0^s (B_{\tau}(s,\tau)\bxq(\tau),\bxq(t))d\tau ds
-\int_0^t\int_0^s (B_{\tau}(s,\tau)\bxq(\tau),\bxq(s))d\tau ds\right|\\
&\leq & \frac{a_1T}{a_0} \left\{\|A^\half\bxq(t)\| \int_0^t\|A^\half\bxq(s)\|ds+ \int_0^t\|A^\half\bxq(s)\|^2ds\right\}.
\end{eqnarray*}
Using the Cauchy-Schwarz inequality and the bounds for $I_1$ and $I_2$, we obtain
\begin{eqnarray*}\label{xiu-2}
\|\xi_{u}(t)\|^2&+&\|A^\half\bxq(t)\|^2\label{xi2-a-n}\leq \|\xi_{u}(0)\|^2+2\int_0^t
\|\eta_{u_{t}}(s)\|\|\xi_{u}(s)\|ds\\
&+&C(T,a_0,a_1)\left(\|A^\half\bxq(t)\|\int_0^t\|A^\half\bxq(s)\|ds+\int_0^t\|A^\half\bxq(s)\|^2ds\right). 
\nonumber   
\end{eqnarray*}
Now, let $|||(\xi_{u},\bxq)(t)|||^2=\|\xi_{u}(t)\|^2+\|A^\half\bxq(t)\|^2$ and 
$$|||(\xi_{u},\bxq)(t^\ast)|||=\max_{0\leq s\leq t}|||(\xi_{u},\bxq)(t)|||,$$
for some $t^\ast \in [0,t]$. Then, at $t=t^\ast$, we note that
\begin{eqnarray*}\label{xiu-2-1a}
|||(\xi_{u},\bxq)(t^\ast)|||&\leq& |||(\xi_{u},\bxq)(0)|||+2\int_0^{t^\ast}
\|\eta_{u_{t}}(s)\|ds\\
&&+C(T,a_0,a_1)\int_0^{t^\ast}|||(\xi_{u},\bxq)(s)||\,ds,\nonumber 
\nonumber   \hspace{-2cm}
\end{eqnarray*}
and therefore,
\begin{eqnarray*}\label{xiu-2-2a}
|||(\xi_{u},\bxq)(t)|||&\leq&|||(\xi_{u},\bxq)(0)|||+2\int_0^{t}\|\eta_{u_{t}}(s)\|\,ds\nonumber\\
&&+C(T,a_0,a_1)\int_0^t|||(\xi_{u},\bxq)(s)|||\,ds.\nonumber 
\nonumber   
\end{eqnarray*}
An application of Gronwall's lemma yields
$$
\|\xi_{u}(t)\|+\|A^\half\bxq(t)\|
\leq C\left(\|\xi_{u}(0)\| +\int_0^t||\eta_{u_{t}}||\,ds\right).
$$
Finally, a use of the triangle inequality with Lemma~\ref{lm3} concludes the proof of 
Theorem~\ref{thm2}. \hfill{$\Box$}

\se

\section{ Error Estimates for a Completely Discrete Scheme}
\se
In this section, we  introduce further notations and formulate 
a completely discrete scheme by applying an implicit finite difference method to 
discretize the time variable of the semidiscrete system (\ref{sw1})-(\ref{sw3}).
Then, we discuss optimal error estimates.

Let $\dt$ $(0<\dt<1)$ be the time step, $k=T/N$ for some positive integer $N$,  and $t_n=n\dt$.
For any function $\phi$ of time, let $\phi^n$ denote $\phi(t_n)$.
We shall use this notation for functions
defined for continuous in time as well as those defined for discrete in time.
Set $$\phi^\nph=\frac{\phi^{n+1}+\phi^n}{2},\qquad
\phi^{n;\qr}=\frac{\phi^{n+1}+2\phi^n+\phi^{n-1}}{4}=\frac{\phi^\nph+\phi^\nmh}{2},$$
and define the following terms for the difference quotients:
$$
\ph \phi^\nph=\frac{\phi^{n+1}-\phi^n}{k},\qquad \pb \phi^\nph=\frac{\phi^\nph-\phi^\nmh}{k},
$$
$$
 \pc \phi^n=\frac{\phi^{n+1}-\phi^{n-1}}{2k}=
\frac{\ph \phi^\nph+\ph \phi^\nmh}{2},
$$
and
$$
\ps \phi^n=\frac{\phi^{n+1}-2\phi^n+\phi^{n-1}}{2k}=\frac{\ph \phi^\nph-\ph \phi^\nmh}{k}.
$$
The discrete-in-time scheme is based on a symmetric difference approximation
around the nodal points, and integral terms are computed by using
the second order quadrature formula
$$
\ii^n(\phi)=k\sum_{j=0}^{n-1}g(t_{j+1/2})\approx\int_0^{t_n}g(s)\,ds,\quad
\mbox{with} \quad t_{j+1/2}=(j+1/2)\dt.
$$
The quadrature error ${\cal E}^n(g)$ is defined by
$$
{\cal E}^n(g)=\ii^{n}(g)-\int_0^{t_n}g(s)\,ds=
\sum_{j=0}^{n-1}\left( kg^{j+\half}-\int_{t_j}^{t_{j+1}}g(s)\,ds \right).
$$
Using Peano's kernel theorem, see \cite{PTW}, it can be written as
$$
{\cal E}^n(g)=\sum_{j=0}^{n-1}\int_{t_j}^{t_{j+1}}\psi(s)D_s^2g(s)\,ds
$$
with
$$
\psi(s)=\left\{\begin{array}{ll} 
(s-t_j)(s-t_{j+\half}),& s\in[t_j,t_{j+\half}],\\
(s-t_{j+1})(s-t_{j+\half}),& s\in[t_{j+\half},t_{j+1}].
\end{array}\right.
$$
Now, let ${\cal B}(t,s;\cdot,\cdot):\bV\times\bV\longrightarrow \R$ be the bilinear form
defined by
$${\cal B}(t,s;{\boldsymbol\phi},{\boldsymbol\chi})=
(B(t,s){\boldsymbol\phi}(s),{\boldsymbol\chi}).$$
Then, the discrete-in-time scheme for the problem
(\ref{w1})-(\ref{w3}) is to seek  
$(U^{n},\Q^{n},\Z^{n})\in W_h\times{\bV}_h\times{\bV}_h$, such that
\begin{eqnarray}
&&\frac{2}{\dt}(\ph U^{\half},w_h)-(\nabla\cdot\Z^{\half},w_h)=(\frac{2}{\dt}u_1,w_h), 
\quad w_h\in W_h,\label{7-1}\\
&&(\Q^{n+\half},\bv_h)+(U^{n+\half},\nabla\cdot\bv_h)=0,\; n\geq 0,\quad \bv_h\in {\bV}_h,\label{7-1a}\\
&&(\Z^{n+\half},{\bz}_h)-(A\Q^{n+\half},{\bz}_h)+
\ii^{n+\half}(\B^{n+\half}(\Q,{\bz}_h))=0,\; n\geq 0,\quad \bz_h\in {\bV}_h,\label{7-1b}\\
&&(\ps U^{n},w_h)-(\nabla\cdot\Z^{n;\qr},w_h)=0,\; n\geq 1,\quad w_h\in W_h,\label{7-1c}
\end{eqnarray}
with given initial data $(U^{0},\Q^{0},\Z^{0})$ in $W_h\times\bV_h\times\bV_h$. Here, 
in (\ref{7-1b}),
$$
\ii ^{n+\half}(\B^{n+\half}(\Q,{\bz}_h))=\frac{1}{2}\left(\ii^{n+1}
(\B^{n+1}(\Q,{\bz}_h))+\ii^n(\B^n(\Q,{\bz}_h)) \right),
$$
where
$$\ii^n(\B^n(\Q,\chi))=\dt\sum_{j=0}^{n-1}(B(t_n,t_{j+1/2})\Q^{j+1/2},\chi).$$
This choice of the time  discretization leads to a second order accuracy in $\dt$.

For ${\boldsymbol\phi}\in\bV_h$, we define a linear functional 
$\ce^{n}({\boldsymbol\phi})$  representing the error in 
 the quadrature formula by
$$
\ce^{n}({\boldsymbol\phi})({\boldsymbol\chi})=\ii^{n}\left(\B^{n}({\boldsymbol\phi},
{\boldsymbol\chi})\right)-\int_0^{t_n}\B(t_n,s;{\boldsymbol\phi},{\boldsymbol\chi})\,ds.
$$
Notice that $\ce^{0}({\boldsymbol\phi})=0$. 
In our analysis, we shall use the following lemma, which can be found in
\cite{PTW}.
\begin{Lemma}\label{lm:8-1}
There exists a positive constant $C$ independent of $h$ and $k$ such that 
the following estimates hold:
$$
k\sum_{n=0}^{m}||\ce^{n+1}({\boldsymbol\phi})||\leq Ck^2\int_0^{t_{m+1}}(||{\boldsymbol\phi}||+
||{\boldsymbol\phi}_{t}||+||{\boldsymbol\phi}_{tt}||)\,ds,
$$
and
$$
k\sum_{n=0}^{m}||\ph\ce^{n+1/2}({\boldsymbol\phi})||\leq Ck^2
\int_0^{t_{m+1}}(||{\boldsymbol\phi}||+||{\boldsymbol\phi}_t||+||{\boldsymbol\phi}_{tt}||)\,ds.
$$
\end{Lemma}
Define $e_U^n:=u^n-U^n$, ${\bf e}_\Q^n:= {\bq^n-\Q^n}$ and ${\bf e}_\Z^n:=\bs^n-\Z^n$. 
From  (\ref{7-1})-(\ref{7-1c})  and (\ref{w1})-(\ref{w3}), we derive the system of equations
\begin{eqnarray}
&&\frac{2}{\dt}(\ph e_U^{\half},w_h)-(\nabla\cdot{\bf e}_\Z^{\half},w_h)=-(2r^0,w_h),\quad w_h\in W_h,\label{8-0-e}\\
&& ({\bf e}_\Q^\nph, \bv_h)+(e_U^\nph,\nabla\cdot\bv_h)=0,\quad \bv_h\in {\bV}_h, \label{8-1-e}\\
&& ({\bf e}_\Z^\nph, {\bz}_h )-(A{\bf e}_\Q^\nph,{\bz}_h)+
  \ii^{n+\half}(\B^{n+\half}({\bf e}_\Q,{\bz}_h))=\ce^{n+1/2}({\bq})({\bz}_h),
  \quad \bz_h\in {\bV}_h,\label{8-2-e}\\
&&(\ps e_U^n, w_h)-(\nabla\cdot{\bf e}_\Z^{n;1/4},w_h)=-(r^n,w_h),\quad w_h\in W_h,\label{8-3-e}
\hspace{-2cm}
\end{eqnarray}
where $\displaystyle r^0=\frac{1}{2}u_{tt}^\half+\frac{1}{\dt}
\left(u_t(0)-\ph u^\half\right)$ and 
\begin{equation}\label{qqw1}
r^n=u_{tt}^{n;1/4}-\ps u^n=\frac{1}{12}\\
\int_{-k}^{k}(|t|-k)\left(3-2(1-|t|/k)^2\right)\frac{\partial^4 u}{\partial t^4}
(t^n+t)dt,\;\;n\geq 1.
\end{equation}
In order to derive {\it a priori} error estimates for the completely discrete scheme, 
we rewrite
\begin{eqnarray*}
e_U^n=(u^n-\tilde u^n_h)-(U^n-\tilde u^n_h)=:\eU^n-\xU^n,\\
         {\bf e}_\Q^n=(\bq^n-\tilde{\bq}^n_h)-(\Q^n-\tilde{\bq}^n_h)=:\eQ^n-\xQ^n,\\
          {\bf e}_\Z^n=(\bs^n-\tilde{\bs}^n_h)-(\Z^n-\tilde{\bs}^n_h)=:\eZ^n-\xZ^n.
\end{eqnarray*}
Since estimates for $\eU$, $\eQ$ and $\eZ$ are known from Lemmas~\ref{lm2} and \ref{lm3}, 
it is sufficient to estimate $\xU$, $\xQ$ and 
$\xZ$. From (\ref{8-0-e})-(\ref{8-3-e}), we obtain the following system
\begin{eqnarray}
&&\frac{2}{\dt}(\ph \xU^{\half},w_h)-(\nabla\cdot\xZ^{\half},w_h)=
\frac{2}{\dt}(\ph \eU^{\half},w_h)+(2r^0,w_h),\label{8-0}\\
&& (\xQ^\nph, \bv_h)+(\xU^\nph,\nabla\cdot\bv_h)=0, \label{8-1}\\
&& (\xZ^\nph, {\bz}_h )-(A\xQ^\nph,{\bz}_h)+
  \ii^{n+\half}(\B^{n+\half}(\xQ,{\bz}_h))=-\ce^{n+1/2}(\tilde{\bq}_h)({\bz}_h),\label{8-2}\\
&&(\ps \xU^n, w_h)-(\nabla\cdot\xZ^{n;1/4},w_h)=
(\ps\eU^n, w_h)+(r^n,w_h).\label{8-3}
\end{eqnarray}
Below, one of the main theorems of this section is proved.
\begin{theorem}\label{thm3} 
Let $\Omega$ be a bounded convex polygonal domain in $\R^2$ and 
let $(u,\bq,\bs)$ be the solution of $(\ref{w1})$-$(\ref{w3})$. 
Further, let $(U^n,\Q^n,\Z^n)\in W_h\times{\bV}_h\times  {\bV}_h$ be the 
solution of $(\ref{7-1})$-$(\ref{7-1c})$. Assume that $U^0=\tilde{u}_h(0)$ 
and $\Q^0=P_h\nabla u_0$.
Then there exists  constants $C>0,$ independent of $h$ and $k$, 
and  $k_0>0,$ such that for $0<k<k_0$ and
$m=0,1,\cdots,N-1$ 
\begin{equation}\label{eu2-d}
\|\ph (u(t_{m+\half})-U^{m+\half})\|
\leq  C(h^2+k^2)\left(\|u_0\|_4+\|u_1\|_3\right),
\end{equation}
and
\begin{equation}\label{eu2-d-b}
\|{\bf q}(t_{m+\half})-\Q^{m+\half}\|+
\|\bs(t_{m+\half})-\Z^{m+\half}\|
\leq C(h^2+k^2)\left(\|u_0\|_4+\|u_1\|_3\right).
\end{equation}
\end{theorem}
Proof. Write (\ref{8-1}) at $n-1/2$ and subtract the resulting one from (\ref{8-1}). Take averages of (\ref{8-2}) 
at two time levels, that is, at $n+1/2$ and $n-1/2.$ Then choose $\bv_h= \xZ^{n;1/4}$ in the modified equation (\ref{8-1}),
${\bz}_h= -\delta_t \xQ^n$ in the modified equation (\ref{8-2}) and $w_h= \delta_t \xU^n$ in (\ref{8-3}). Then add to obtain
\begin{eqnarray}\label{8-main}
\frac{1}{2} \pb \Big(||\ph \xU^{n+\half}||^2 &+& ||\ah \xQ^{n+\half}||^2\Big) = (\ps\eU^n + r^n, \delta_t\xU^n)+ \ce^{n;1/4}(\tilde{\bq}_h)(\delta_t \xQ^n) \nonumber\\
&+&  \frac{1}{2} \Big(\ii^{n+\half}(\B^{n+\half}(\xQ,\delta_t \xQ^n ))+ \ii^{n-\half}(\B^{n-\half}(\xQ,\delta_t \xQ^n )) \Big)\\
&=& I_1^n+ I_2^n + I_3^n. \nonumber
\end{eqnarray}
With $\Psi^n=(\xU^n,\xQ^n )$, define
$$|||\Psi^{n+\half}|||^2=||\ph \xU^{n+\half}||^2+||\ah \xQ^{n+\half}||^2.$$
Multiply (\ref{8-main}) by $2 k$ and sum from $n=2$ to $m$ to obtain
\begin{equation}\label{8-main-1}
|||\Psi^{m+1/2}|||^2
\leq |||\Psi^{3/2}|||^2
+2\dt\left|\sum_{n=2}^m(I_1^n+I_2^n+I_3^n)\right|.
\end{equation}
Define for some $m^{\star}$ with $ 0 \le m^{\star}\leq m,$
$$|||\Psi^{m^{\star}+1/2}||| =\max_{0\leq n\leq m}|||\Psi^\nph|||.$$
A use of the Cauchy-Schwarz inequality  yields
\begin{eqnarray*}
2\dt\left|\sum_{n=2}^m I_1^n \right|
&\leq & \sum_{n=2}^m\left(\|\partial^2_t \eta_u^n\|+\|r^n\|\right)\, 
\left(\|\ph \xU^\nph\|+\|\ph \xU^\nmh\| \right)\\
&\leq & 2\dt \sum_{n=2}^{m}\left(\|\partial^2_t \eta_u^n\|+\|r^n\|\right)|||\Psi^{m^{\star}+1/2}|||.
\end{eqnarray*} 
Setting 
$$ \tilde{B}^{n+\half}_{j+\half}=\frac{1}{2}\left(B (t_{n+1},t_\jph)+ B (t_{n},t_\jph) \right),$$
we now estimate one of the term of $I_3^n.$  Note that
\begin{equation*}
\begin{split}
\ii^{n+\half}&(\B^{n+\half}(\xQ,\pc\xQ^n))\\
=&\frac{\dt}{2}\left[\sum_{j=0}^{n} \B(t_{n+1},t_{j+1/2};\xQ^{j+1/2},\pc \xQ^n)
+\sum_{j=0}^{n-1} \B(t_n,t_{j+1/2};\xQ^{j+1/2},\pc \xQ^n)\right]\\
=&\frac{1}{2}  \Big(B(t_{n+1},t_\nph)\xQ^\nph, (\xQ^\nph-\xQ^\nmh)\Big)
+\sum_{j=0}^{n-1} \Big(\Bnj \xQ^\jph,(\xQ^\nph-\xQ^\nmh)\Big)\\
=&I_{31}^n+I_{32}^n.
\end{split}
\end{equation*}
For the second term $I_{32}^n$, we use the fact that
\begin{equation}\label{eq:formula}
H^\nph\pb\xQ^\nph=\pb(H^\nph\xQ^\nph)-\pb(H^\nph)\xQ^\nmh,
\end{equation}
to obtain 
\begin{eqnarray*}
I_{32}^n&=& k\pb\left(\sum_{j=0}^{n-1} \Big(\Bnj \xQ^\nph, \xQ^\jph\Big)\right)
+\Big({\tilde B}^{n-\half}_{n-\half}\xQ^\nmh,\xQ^\nmh\Big)\\
&&-k\sum_{j=0}^{n-1}\Big(\pb\left(\Bnj\right)\xQ^\jph, \xQ^\nmh\Big),
\end{eqnarray*}
and hence after summing up from $n=2$ to $m$ and multiplying by $k$
\begin{eqnarray*}
k\left|\sum_{n=2}^{m} I_{32}^n \right|&\leq&
k^2\left|\sum_{j=0}^{m-1} \Big( {\tilde B}^{m+\half}_\jph \xQ^\jph, \xQ^{m+\half}\Big)
-\Big( {\tilde B}^{3/2}_{1/2}\xQ^{3/2},\xQ^{1/2}\Big)\right|
+k\left|\sum_{n=2}^{m}\Big({\tilde B}^{n-\half}_{n-\half}\xQ^\nmh,\xQ^\nmh\Big)\right|\\
&&+k^2\left|\sum_{n=2}^{m}\sum_{j=0}^{n-1}\Big( \pb\left(\Bnj\right)\xQ^\jph,\xQ^\nmh\Big)\right| \\
&\leq& a_1 k^2 \|\xQ^{m+\half}\|\sum_{j=0}^{m-1}\|\xQ^\jph\|+a_1 k^2\| \xQ^{3/2}\|\,\|\xQ^{1/2}\|
+a_1 (1+T) k\sum_{j=0}^{m-1} \|\xQ^\jph\|^2\\
&\leq& C(a_0,a_1,T)k\left(\sum_{n=0}^{m-1}\|A^{\half}\xQ^\nph\|\right)|||\Psi^{\msph}|||.
\end{eqnarray*}
All together,  we obtain 
\begin{eqnarray*}
k\left|\sum_{n=2}^{m} \ii^{n+\half} (\B^{n+\half}(\xQ,\pc\xQ^n))  \right| 
&\leq& \frac{a_1}{a_0^2} k \|\xQ^{m+\half}\|^2 +
C(a_0,a_1,T)k\left(\sum_{n=0}^{m-1}\|A^{\half}\xQ^\nph\|\right)|||\Psi^{\msph}|||.
\end{eqnarray*}
We can now estimate $\ii^{n-\half}(\B^{n-\half}(\xQ,\pc\xQ^n))$ in a similar way, but without 
having the term $\|A^{\half}\xQ^{m+\half}\|^2$ on the right hand side and thus, we arrive at 
$$
k\left|\sum_{n=2}^{m}I_3^n\right|\leq 
\frac{a_1k}{a_0^2}\|A^{\half}\Q^{m+\half}\|^2
+ C(a_1,a_0,T)k\left(\sum_{n=0}^{m-1}\|A^{\half}\Q^\nph\|\right)|||\Psi^{\msph}|||.
$$
In a similar way, using (\ref{eq:formula}),   the term  $I_2^n$ can be estimated as follows
$$
k\left|\sum_{n=2}^{m}I_2^n\right|
\leq\left|\ce^{n;1/4}(\tilde{\bq}_h)(\xQ^{m+1/2})-\ce^{1;1/4}(\tilde{\bq}_h)(\xQ^{3/2})-
k\sum_{n=2}^{m}\pb (\ce^{n;1/4}(\tilde{\bq}_h))(\xQ^\nmh)\right|.$$
Notice that, since $\ce^0=0$, it follows that
$$
\ce^{n}(\tilde{\bq}_h)(\xQ^n)=k\sum_{j=0}^{n-1}\ph \ce^{j+\half}(\tilde{\bq}_h)(\xQ^n), 
$$
and hence, we obtain 
$$
k\left|\sum_{n=2}^{m}I_2^n\right|\leq C(a_0,a_1)k
\left(\sum_{n=0}^{m}\|\ph \ce^{n+\half}(\tilde{\bq}_h)\| \right)\;|||\Psi^{\msph}|||.
$$
It remains now to bound the  term $|||\Psi^{3/2}|||$ on the right hand side of (\ref{8-main-1}).
Take (\ref{8-main}) at $n=1$ to obtain
\begin{eqnarray}\label{8-3b-nc}
|||\Psi^{3/2}|||^2&\leq&|||\Psi^{1/2}|||^2+2k\left|(\partial_t^2\eta_U^1+ r^1,\pc \xU^1)+
\frac{1}{2}\left(\ii^{3/2}(\B^{3/2}(\xQ,\pc\xQ^1))+\ii^{1/2}(\B^{1/2}(\xQ,\pc\xQ^1)\right)\right.\nonumber\\
&+&\left. \ce^{1;1/4}(\tilde{\bq}_h)(\pc\xQ^1)\right|\nonumber\\
&\leq&|||\Psi^{1/2}|||^2+C(a_0,a_1)k\Big(||\partial_t^2\eta_U^1|| +\|r^1\|
+||A^\half\xQ^{1/2}||\\
&&\left.+||A^\half\xQ^{3/2}||+||\ph \ce^{1/2}(\tilde{\bq}_h)||+||\ph \ce^{3/2}(\tilde{\bq}_h)||\right)
|||\Psi^{\msph}|||.\nonumber
\end{eqnarray} 
Substitution of  estimates involving $I_1^n,\cdots,I_3^n$ and (\ref{8-3b-nc})  in (\ref{8-main-1}) yields 
\begin{eqnarray}\label{8-3c}
\Big(1-\frac{a_1}{a_0^2}k\Big)\,|||\Psi^{m+\half}|||^2 &\leq&  |||\Psi^{1/2}|||^2
+Ck\Big\{\sum_{n=1}^{m}\Big(\|\partial_t^2\eta_U^n\|+\|r^n\|\Big)\nonumber\\
&+&\sum_{n=0}^{m}\|\ph \ce^{n+\half}(\tilde{\bq}_h)\|+\sum_{n=0}^{m-1}|||\Psi^\nph|||\Big\}
|||\Psi^{\msph}|||,
\end{eqnarray}
Choose $k_0>0$ such that for $0 < k \leq k_0,$ $ (1-\frac{a_1}{a_0^2 }k )>0.$ Then 
replace $m$ by $m^{\star}$ in (\ref{8-3c}) and obtain after cancellation of $|||\Psi^{\msph}|||$  from the both sides
\begin{eqnarray*}
|||\Psi^{m+\half}|||&\leq& |||\Psi^{\msph}||| \leq   C\Big\{ |||\Psi^{1/2}|||
+ k \sum_{n=1}^{m}\Big(\|\partial_t^2\eta_U^n\|+\|r^n\|\Big) \\
&+& \sum_{n=0}^{m}\|\ph \ce^{n+\half}(\tilde{\bq}_h)\|+\sum_{n=0}^{m-1}|||\Psi^\nph|||\Big\}.\nonumber
\end{eqnarray*} 
Then an application of the discrete Gronwall's lemma yields
\begin{eqnarray}
|||\Psi^{m+\half}|||\leq C\left\{
|||\Psi^\half|||+k\sum_{n=1}^{m}\big(||\ps\eU^n||+||r^n||\big)
+k\sum_{n=0}^{m}||\ph \ce^{n+1/2}(\tilde{\bq}_h)||\right\}.
\end{eqnarray}
To estimate $\Psi^\half$ on the right hand side 
of this inequality, we first note that, 
$\xU^0=0$ since $U^0=\tilde{u}_h(0)$, and hence, 
$\displaystyle\xU^\half=\frac{k}{2}\ph \xU^\half$.
Now,  choose $w_h=\xU^\half$ 
in (\ref{8-0}), $w_h=\xZ^\half$ in (\ref{8-1}), and $w_h=-\xQ^\half$ in (\ref{8-2}).
Adding the resulting equations, and taking into account that $\ii^0(\phi)=0$ 
and $\ce^0(\phi)=0$, we arrive at
\begin{eqnarray*}
|||\Psi^\half|||^2
&\leq&
\left|\left(\ph\eU^\half,\ph \xU^\half\right)\right|+k\left|\left(r^0,\ph \xU^\half\right)\right|
+\frac{1}{2}\left|\ii^1(\B^1(\xQ,\xQ^\half)\right|
+\frac{1}{2}\left|\ce^1(\tilde{\bq}_h)(\xQ^\half)\right|\\
&\leq& C(a_0)
\left(||\ph\eU^\half||+k||r^0||+||\ce^1(\tilde{\bq}_h)||\right)|||\Psi^\half||| +
\frac{a_1k}{2a_0}||A^\half\xQ^\half||^2 \\
&\leq& C(a_0)
\left(||\ph\eU^\half||+k||r^0||+||\ce^1(\tilde{\bq}_h)||\right)|||\Psi^\half||| +
\frac{a_1k}{2a_0}|||\Psi^\half|||^2
\end{eqnarray*}
For $ 0<k\leq k_0,$ $\big(1-(a_1k)/(2a_0)\big) >0 $ and hence, we obtain
$$
|||\Psi^\half|||\leq C\left\{||\ph\eU^\half||+k||r^0||+||\ce^1(\tilde{\bq}_h)||\right\}.
$$
Thus,
\begin{equation}\label{8-5}
|||\Psi^{m+\half}|||\leq C\left\{
||\ph\eU^\half||+k\sum_{n=1}^{m}||\ps\eU^n||+k\sum_{n=0}^{m}||r^n||
+k\sum_{n=0}^{m}||\ph \ce^{n+1/2}(\tilde{\bq}_h)||
\right\}.
\end{equation}
To estimate the first two terms on the right hand side of (\ref{8-5}), it is observed that
\begin{eqnarray}\label{8-5a}
||\ph\eU^\half||\leq\frac{1}{k}\int_0^{\dt}||\eta_{\U t}(s)||\,ds,
\end{eqnarray}
and a use of Taylor series expansions yields
\begin{eqnarray}\label{8-5b-n}
k\sum_{n=1}^{m}||\ps\eU^n||&\leq&\frac{1}{k}\sum_{n=1}^{m}\left\{
\int_{t_n}^{t_{n+1}}(t_{n+1}-s)||\eta_{\U tt}(s)||\,ds+
\int_{t_{n-1}}^{t_n}(s-t_{n-1})||\eta_{\U tt}(s)||\,ds
\right\}\nonumber\\
&\leq&2\int_{0}^{t_{m+1}}||\eta_{\U tt}(s)||\,ds\leq C(T)h^2(\|u_0\|_4+\|u_1\|_3).
\end{eqnarray}
Further, from (\ref{qqw1}) it follows that
$$||r^n||\leq Ck\int_{t_{n-1}}^{t_{n+1}}\|D_t^4u(s)\|\,ds,\quad n\geq 1,$$
and 
$$||r^0||\leq Ck||u_{ttt}||_{L^\infty(0,\dt/2;L^2(\Omega))}
\leq Ck\int_{0}^{t_{m+1}}
\left( \|D_t^3u(s)\|+\|D_t^4u(s)\|\right)\,ds.$$
Hence,
\begin{equation}\label{s2}
k\sum_{n=0}^{m}||r^n||\leq Ck^2\int_{0}^{t_{m+1}}
\left( \|D_t^3u(s)\|+\|D_t^4u(s)\|\right)\,ds
\leq C(T)k^2(\|u_0\|_4+\|u_1\|_3).
\end{equation}
For the last term in (\ref{8-5}), a use of Lemma~\ref{lm:8-1} with the triangle
inequality yields
$$
k\sum_{n=0}^{m}||\ph\ce^{n+1/2}(\tilde{\bq}_h)||\leq Ck^2
\sum_{j=0}^2\int_{0}^{t_{m+1}} \left(\|D_t^j \bq(s)\|+\|D_t^j \eq(s)\|\right)\,ds,
$$
and hence,  the estimates  in Lemmas~\ref{lm1} and \ref{lm2} show that
\begin{equation}\label{s3}
k\sum_{n=0}^{m}||\ph\ce^{n+1/2}(\tilde{\bq}_h)||\leq 
C(T)k^2(\|u_0\|_3+\|u_1\|_2)+C(T)k^2h(\|u_0\|_4+\|u_1\|_3).
\end{equation}
Substitute (\ref{8-5a})-(\ref{s3}) in (\ref{8-5}) and use the triangle inequality 
with the results in Lemmas~\ref{lm2} and \ref{lm3}, 
to obtain the error estimates involving $u$ and $\bq$ in (\ref{eu2-d})-(\ref{eu2-d-b}). 
Now to complete the remaining estimate in (\ref{eu2-d-b}), we need to estimate 
$||\xZ^{m+\half}||$. To do so, choose
${\bz}_h=\xZ^{m+\half}$ in (\ref{8-2}), and conclude that
$$||\xZ^{m+\half}||\leq C(a_1,T)\left(\max_{0\leq n\leq m}||\xQ^{n+\half}|| + 
||\ce^{m+1/2}(\tilde{\bq}_h)||\right).$$
Finally, a use of the triangle inequality with Lemmas~\ref{lm2} and \ref{lm3}
completes the rest of the proof. 
\hfill$\Box$


Below, we again recall a variant of Baker's nonstandard energy formulation  
to prove  $\ell^{\infty}(L^2)$-estimate for the error $\|e_U^n\|$ under reduced regularity conditions on the initial data.

Now introduce the following notations for proving the next theorem.
Define 
$$
\hat{\phi}^0=0,\qquad \hat{\phi}^n=k\sum_{j=0}^{n-1}\phi^{\jph}.
$$
Then, 
$$\ph \hat{\phi}^\nph={\phi}^\nph,$$ and 
$$k\sum_{j=1}^{n}\phi^{j;1/4}=\hat{\phi}^{\nph}-\frac{k}{2}\phi^{\half}.$$
Notice that 
\begin{eqnarray}\label{it1}
k\sum_{n=0}^m\psi^{n+1}\phi^\nph&=&k\sum_{n=0}^m\psi^{n+1}\ph \hat{\phi}^\nph\\
&=&\sum_{n=0}^m(\psi^{n+1} \hat{\phi}^{n+1}-\psi^{n}\hat{\phi}^{n})
-\sum_{n=0}^m(\psi^{n+1}-\psi^{n})\hat{\phi}^{n}\nonumber\\
&=&\psi^{m+1}\hat{\phi}^{m+1}-k\sum_{n=0}^m\ph\psi^\nph \hat{\phi}^{n}.\nonumber
\end{eqnarray}
and similarly, 
\begin{equation}\label{it2}
k\sum_{n=0}^m\psi^\nph\phi^\nph=\psi^{m+\half}\hat{\phi}^{m+1}-
k\sum_{n=0}^m\pb\psi^\nph \hat{\phi}^{n}.
\end{equation}
Multiplying (\ref{7-1c}) by $k$, summing over $n$, 
and taking into account (\ref{7-1}), we obtain the new equation
\begin{equation}\label{7-1cn}
(\ph U^{n+\half},w_h)-(\nabla\cdot\hat{\Z}^{n+\half},w_h)=(u_1,w_h).
\end{equation}
The key idea here, which differs from Baker's approach, is that we compare
the above equation with
\begin{equation}\label{w3n}
(u_{t},w)-(\nabla\cdot\bar{\bs},w) = (u_1,w)
\end{equation}
which is derived by integrating (\ref{w3}) with respect to time, 
where $\bar{\bs}(t)=\int_0^t\bs(s)ds$.
By taking the average of (\ref{w3n}) at  $t^{n+1}$ and $t^{n}$ and using (\ref{7-1cn})
we arrive at
$$
(\ph e_U^{\nph},w)-(\nabla\cdot\hat{\bs}^{n+\half}-\nabla\cdot\hat{\Z}^{n+\half},w)=-(r^n_1,w)
-(\nabla\cdot\hat{\bs}^{n+\half}-\nabla\cdot\bar{\bs}^\nph,w),
$$
with $r^n_1=u_t^{\nph}-\ph u^{\nph}$. Thus,
$$
(\ph e_U^{\nph},w)-(\nabla\cdot{\hat{\bf e}_\Z}^{n+\half},w) = -(r^n_1,w)
-{\cal E_I}^\nph(\nabla\cdot\bs)(w_h),
$$
where
$$
{\cal E}^n_{\cal I}(\phi)(\chi)=(\hat{\phi}^{n}-\bar{\phi}^n,\chi)
=\dt\sum_{j=0}^{n-1}(\phi^{j+1/2},\chi)-\int_0^{t_n}(\phi(s),\chi)\,ds.
$$
Using the definitions of $\tilde{u}_h$ and $\tilde{\bs}_h$, it follows that 
$$
(\ph\xU^\nph, w_h)-(\nabla\cdot\xZh^\nph,w_h)=
(\ph\eU^\nph, w_h)+(r_1^n,w_h)+{\cal E}^\nph_{\cal I}(\nabla\cdot\bs)(w_h).
$$
To complete the system of error equations, we multiply (\ref{8-2}) by $k$, sum over $n$, 
and take the average of the resulting equations. Including (\ref{8-2}), we end up
with following system 
\begin{eqnarray}
&& (\xQ^\nph, \bv_h)+(\xU^\nph,\nabla\cdot\bv_h)=0, \label{8-1b}\\
&& (\xZh^\nph, {\bz}_h)-(A\xQh^\nph,{\bz}_h)+\hat{\ii}^{n+\half}(\xQ,{\bz}_h)=
-\hce^{n+1/2}(\tilde{\bq}_h)({\bz}_h),\label{8-2b}\\
&&(\ph\xU^\nph, w_h)-(\nabla\cdot\xZh^\nph,w_h)=
(\ph\eU^\nph, w_h)+(r_1^n,w_h)+{\cal E}^\nph_{\cal I}(\nabla\cdot\bs)(w_h),\label{8-3b}
\end{eqnarray}
where 
$$
\hat{\ii}^{n+\half}(\xQ,{\bz}_h)=\frac{1}{2}\left[
k\sum_{j=0}^{n}\ii^{j+\half}(\B^{j+\half}(\xQ,{\bz}_h))+
k\sum_{j=0}^{n-1}\ii^{j+\half}(\B^{j+\half}(\xQ,{\bz}_h))\right],
$$
and
$$
\hce^{n+1/2}(\tilde{\bq}_h)({\bz}_h)=\frac{1}{2}\left[  
k\sum_{j=0}^{n}\ce^{j+1/2}(\tilde{\bq}_h)({\bz}_h)+
k\sum_{j=0}^{n-1}\ce^{j+1/2}(\tilde{\bq}_h)({\bz}_h)\right].
$$

\begin{theorem}\label{thm4}
Let $\Omega$ be a bounded convex polygonal domain in $\R^2$, and 
let $(u,\bq,\bs)$ be the solution of $(\ref{w1})$-$(\ref{w3})$. 
Further, let $(U^n,\Q^n,\Z^n)\in W_h\times{\bV}_h\times  {\bV}_h$ be the 
solution of $(\ref{7-1})$-$(\ref{7-1c})$. With $U^0=P_hu_0$ and  
$\Q^0=P_h\nabla u_0$, there exists a positive constant $C$, independent of $h$ and $k$,
such that for small $k$ with $k=O(h)$, the following estimate holds:
\begin{equation}\label{6-2-bb}
\|u(t_{m+1})-U^{m+1}\|
\leq  C (h^2+k^2)\left(\|u_0\|_3+\|u_1\|_2\right),\quad m=0,1,\cdots,N-1.
\end{equation}
\end{theorem}
Proof.  Choose $w_h=\xZh^\nph$ 
in (\ref{8-1b}), ${\bz}_h=-\xQ^\nph$ in (\ref{8-2b}), and 
$v_h=\xU^\nph$ in (\ref{8-3b}).
Adding the resulting equations, we find that
\begin{eqnarray*} \label{8-4b}
(\ph\xU^\nph,\xU^\nph)&+&(A\xQh^\nph,\xQ^\nph) = (\ph\eU^\nph,\xU^\nph)+(r_1^n,\xU^\nph) \nonumber\\
&-& \hat{\ii}^{n+\half}(\xQ,\xQ^\nph)
-\hce^{n+1/2}(\tilde{\bq}_h)(\xQ^\nph) 
+{\cal E}^\nph_{\cal I}(\nabla\cdot\bs)(\xU^\nph).
\end{eqnarray*}
Multiply both sides by $2k$ and sum from $n=0$ to $m$. Then,  use that 
$\xQh^0=0$, $\xQ^\nph=\ph \xQh^\nph,$ and   
$$
(\ph\xU^\nph, \xU^\nph)=\frac{1}{2k}
\left(||\xU^{n+1}||^2-||\xU^{n}||^2\right),
$$ 
to arrive at
\begin{eqnarray}\label{8-5b}
||\xU^{m+1}||^2 &+& \|A^\half\xQh^{m+1}\|^2 =||\xi_U^0||^2
-2k\sum_{n=0}^{m}\hat{\ii}^{n+\half}(\xQ,\xQ^\nph)-2k\sum_{n=0}^{m}\hce^{n+1/2}(\tilde{\bq}_h)(\xQ^\nph)   \nonumber\\
&&+2k\sum_{n=0}^{m}(\ph\eU^\nph,\xU^\nph) +2k\sum_{n=0}^{m}(r_1^n,\xU^\nph)
+2k\sum_{n=0}^{m}{\cal E}^\nph_{\cal I}(\nabla\cdot\bs)(\xU^\nph)\nonumber\\
&=&||\xi_U^0||^2+I_1^m+I_2^m+I_3^m+I_4^m+I_5^m.
\end{eqnarray}
For convenience, we use the following notations: $B_r^s=B(t_{s},t_r)$,
$\bar\partial_i\phi_{i+\half}=(\phi_{i+\half}-\phi_{i-\half})/k$,
and denote by $\bf\psi\cdot \bf\phi$ the $L^2$ inner product of $\bf\psi$ and $\bf\phi$.
Further, let $|||(\xU^n,\xQh^n)|||^2=||\xU^n||^2+||A^\half\xQh^n||^2$ and  for some $m^{\star}\in [0;m+1]$, define
 $$|||(\xU^{m^{\star}},\xQh^{m^{\star}})|||=\max_{0\leq n\leq m+1}|||(\xU^n,\xQh^n)|||.$$
To estimate the first term $I_1^m$, we observe that by (\ref{it2})
\begin{eqnarray*}
\ii^{j+1}(B^{j+1}(\xQ,\xQ^\nph))
&=& k\sum_{i=0}^jB^{j+1}_\iph\xQ^\iph\cdot\xQ^\nph\\
&=&B^{j+1}_\jph\xQh^{j+1}\cdot\xQ^\nph-k\sum_{i=0}^j
\left(\bar\partial_i(B^{j+1}_\iph)\right)
\xQh^i\cdot\xQ^\nph.
\end{eqnarray*}
Setting 
$$\Theta^{n+1}=\sum_{j=0}^nB^{j+1}_\jph\xQh^{j+1}\quad\mbox{ and } \quad
\Upsilon^{n+1}=k\sum_{j=0}^n\sum_{i=0}^j\left(\bar\partial_i(B^{j+1}_\iph)\right)\xQh^i,
$$
we now obtain 
\begin{equation}\label{sd-0}
k^2 \sum_{n=0}^m\sum_{j=0}^n\ii^{j+1}(B^{j+1}(\xQ,\xQ^\nph))=
k^2\sum_{n=0}^m\Theta^{n+1}\cdot\xQ^\nph-k^2\sum_{n=0}^m\Upsilon^{n+1}\cdot\xQ^\nph.
\end{equation}
Next, we estimate the terms $\Theta^{n+1}$ and $\Upsilon^{n+1}$. Using (\ref{it1}) we
have 
$$
k^2\sum_{n=0}^m\Theta^{n+1}\cdot\xQ^\nph
=k\Theta^{m+1}\cdot\xQh^{m+1}-k^2\sum_{n=0}^m\ph \Theta^\nph\cdot\xQh^{n},
$$
and hence, using the definition of $\Theta^n$, it follows that
\begin{equation}\label{sd-1}
k^2\sum_{n=0}^m\Theta^{n+1}\cdot\xQ^\nph
=k\sum_{n=0}^m B^{n+1}_\nph \xQh^{n+1}\cdot\xQh^{m+1}-k\sum_{n=0}^mB^{n+1}_\nph
\xQh^{n+1}\cdot\xQh^{n}.
\end{equation}
Similarly,
\begin{eqnarray*}\label{sd-2}
k^2\sum_{n=0}^m\Upsilon^{n+1}\cdot\xQ^\nph
&=&k\Upsilon^{m+1}\cdot\xQh^{m+1}-k^2\sum_{n=0}^m\ph\Upsilon^{n+\half}\cdot\xQh^{n}\\
&=&k^2\left(\sum_{j=0}^m\sum_{i=0}^j\left(\bar\partial_i(B^{j+1}_\iph)\right)
\xQh^i\right)\cdot\xQh^{m+1}\nonumber\\
&&-k^2\sum_{n=0}^m\left(\sum_{i=0}^n\left(\bar\partial_i(B^{n+1}_\iph)\right)\xQh^i
 \right)\cdot\xQh^{n}.\nonumber
\end{eqnarray*}
On substitution of (\ref{sd-1}) and (\ref{sd-2}) in (\ref{sd-0}), a use of Cauchy-Schwarz
inequality with $||D_{t,s}B(t,s)||\leq a_1$, yields
$$
\left|k^2 \sum_{n=0}^m\sum_{j=0}^n\ii^{j+1}(B^{j+1}(\xQ,\xQ^\nph))\right|
\leq \frac{a_1k}{a_0} ||A^\half\xQh^{m+1}||^2+
C(a_0,a_1,T)\left(k\sum_{n=0}^m ||A^\half\xQh^{n}||\right)|||(\xU^{m^{\star}},\xQh^{m^{\star}})|||.
$$
Since, similar bounds can be obtained for other terms in 
$k\sum_{n=0}^{m}\hat{\ii}^{n+\half}(\xQ,\xQ^\nph)$, we finally conclude
that
$$
|I_1^m|\leq \frac{a_1k}{2a_0} ||A^\half\xQh^{m+1}||+
C(a_0,a_1,T)\left(k\sum_{n=0}^m ||A^\half\xQh^{n}||\right)|||(\xU^{m^{\star}},\xQh^{m^{\star}})|||.
$$
Now, with $\Lambda^{n+1}=\sum_{j=0}^{n}\ce^{j+1}(\tilde{\bq}_h)$, we observe that
\begin{eqnarray*}
k^2\sum_{n=0}^{m}\sum_{j=0}^{n}\ce^{j+1}(\tilde{\bq}_h)(\xQ^\nph)
&=&k^2\sum_{n=0}^{m}\Lambda^{n+1}\ph\xQh^\nph\\
&=&k\Lambda^{m+1}(\xQh^{m+1})-k\sum_{n=0}^m(\Lambda^{n+1}-\Lambda^{n})(\xQh^{n})\\
&=&k\sum_{j=0}^{m}\ce^{j+1}(\tilde{\bq}_h)(\xQh^{m+1})
-k\sum_{n=0}^{m}\ce^{n+1}(\tilde{\bq}_h)(\xQh^n).
\end{eqnarray*}
Since, the terms in $k\sum_{n=0}^{m}\hce^{n+1/2}(\tilde{\bq}_h)(\xQ^\nph)$
have a similar form, we deduce that 
$$
|I_2^m|\leq C(a_0)k\sum_{n=0}^{m}\left(||\ce^{n+1}(\tilde{\bq}_h)||\right)
|||(\xU^{m^{\star}},\xQh^{m^{\star}})|||.
$$
For $I_3^m$ and $I_4^m$, we have 
$$
|I_3^m+I_4^m|\leq C(a_0)k \sum_{n=0}^{m}\left(\left|\left|\ph\eU^\nph\right|\right|+||r_1^n||\right)
|||(\xU^{m^{\star}},\xQh^{m^{\star}})|||.
$$
To estimate the last term $I_5^m$, we first notice that
$(\nabla\cdot\bs,\xU^\nph)=(\nabla\cdot\tilde{\bs}_h,\xU^\nph)$ 
since $\xU^\nph\in W_h$.
Then  by (\ref{8-1b}) 
$$(\nabla\cdot\bs,\xU^\nph)=-(\xQ^\nph,\tilde{\bs}_h).$$
Thus, it follows that
\begin{eqnarray*}
-k\sum_{n=0}^m{\cal E}^{n+1}_{\cal I}(\nabla\cdot\bs)(\xU^\nph)
&=&k\sum_{n=0}^m{\cal E}^{n+1}_{\cal I}(\tilde{\bs}_h)(\xQ^\nph)\\
&=&{\cal E}^{m+1}_{\cal I}(\tilde{\bs}_h)(\xQh^{m+1})-
k\sum_{n=0}^m\ph{\cal E}^\nph_{\cal I}(\tilde{\bs}_h)(\xQh^{n}).
\end{eqnarray*}
Similarly, we have
$$-k\sum_{n=0}^m{\cal E}^{n}_{\cal I}(\nabla\cdot\bs)(\xU^\nph)=
{\cal E}^{m}_{\cal I}(\tilde{\bs}_h)(\xQh^{m+1})-
k\sum_{n=1}^m\ph{\cal E}^\nmh_{\cal I}(\tilde{\bs}_h)(\xQh^{n}),
$$
which yields 
$$|I_5^m|\leq C(a_0)\left(\|{\cal E}^{m+\half}_{\cal I}(\tilde{\bs}_h)\|+
k\sum_{n=0}^m\|\ph{\cal E}^\nph_{\cal I}(\tilde{\bs}_h)\right)|||(\xU^{m^{\star}},\xQh^{m^{\star}})|||.
$$
On substituting the above estimates in (\ref{8-5b}) and following steps  
in previous theorems, we arrive at
\begin{eqnarray*}
\Big(1-(a_1/2a_0)k\Big)\,|||(\xU^{m+1},\xQh^{m+1})||| &\leq &
Ck \sum_{n=0}^{m}\left(
\left|\left|\ph\eU^\nph\right|\right|+||r_1^n||+||\ce^{n+1}(\tilde{\bq}_h)||+
\|\ph{\cal E}^\nph_{\cal I}(\tilde{\bs}_h)\|\right.\\
&&+\left.\|\ph{\cal E}^\nph_{\cal I}(\tilde{\bs}_h)\|  \right)
+C\|{\cal E_I}^{m+\half}(\tilde{\bs}_h)\|+\|\xi_U^0\| \\
 &\leq & Ck \sum_{n=0}^{m}\left(
\left|\left|\ph\eU^\nph\right|\right|+||r_1^n||+||\ce^{n+1}(\tilde{\bq}_h)||
+\|\ph{\cal E}^\nph_{\cal I}(\tilde{\bs}_h)\|\right.\\
&&+\left.|||(\xU^{n},\xQh^{n})|||\right)
+C\|{\cal E}^{m+\half}_{\cal I}(\tilde{\bs}_h)\|+\|\xi_U^0\|.
\end{eqnarray*}
Since for $0 < k \leq k_0,$ $\big(1-(a_1/2a_0)k\big)$ can be made positive, 
an application of the discrete Gronwall's lemma yields
\begin{eqnarray}\label{ss1}
||\xU^{m+1}||+||A^\half\xQh^{m+1}||&\leq&
Ck \sum_{n=0}^{m}\left(
\left|\left|\ph\eU^\nph\right|\right|+\|r_1^n\|+\|\ce^{n+1}(\tilde{\bq}_h)\|
+\|\ph{\cal E}^\nph_{\cal I}(\tilde{\bs}_h)\|\right)\\
&&+C\left(\|{\cal E}^{m+\half}_{\cal I}(\tilde{\bs}_h)\|+\|\xi_U^0\|\right).\nonumber
\end{eqnarray}
The first two terms on the right hand side can be bounded as follows: 
$$
k\sum_{n=0}^{m}\left|\left|\ph\eU^\nph\right|\right|\leq
\int_0^{t_{m+1}}\left|\left|\eta_{\U t}(s)\right|\right|\,ds,
$$ 
and
$$
k\sum_{n=0}^{m}||r^n||\leq Ck^2\int_0^{t_{m+1}}\|D_t^3u(s)\|\,ds.
$$
For the third term, we note that 
$\tilde{\bq}_h=-(\bq-\tilde{\bq}_h)+\bq=-\eQ+\bq$, and hence
\begin{equation}\label{ss2}
k \sum_{n=0}^{m}\|\ce^{n+1}(\tilde{\bq}_h)\|\leq
k \sum_{n=0}^{m}\|\ce^{n+1}(\eQ)\|+k\sum_{n=0}^{m}\|\ce^{n+1}(\bq)\|.
\end{equation}
For the last term on the right hand side of (\ref{ss2}), use Lemmas~\ref{lm:8-1}
and \ref{lm1} to obtain
\begin{eqnarray}\label{ss3}
k\sum_{n=0}^{m}\|\ce^{n+1}(\bq)\|&\leq&
Ck^2\int_0^{t_{m+1}}(||\bq||+
||\bq_{t}||+||\bq_{tt}||)\,ds\nonumber \\
&\leq& Ck^2(\|u_0\|_3+\|u_1\|_2).
\end{eqnarray}
For the first term on the right hand side of (\ref{ss2}), we note that 
\begin{eqnarray}\label{ss4}
\ce^{n+1}(\eQ)({\boldsymbol\chi})&=&\ii^{n+1}\left(\B^{n+1}({\eQ},
{\boldsymbol\chi})\right)-\int_0^{t_{n+1}}\B(t_{n+1},s;\eQ,{\boldsymbol\chi})\,ds\nonumber \\
&=&\sum_{j=0}^n\left[k\left(B(t_{n+1},t_{j+\half})\eQ^{j+\half},{\boldsymbol\chi}\right)
-\int_{t_j}^{t_{j+1}}\B(t_{n+1},s;\eQ,{\boldsymbol\chi})\,ds\right]\nonumber \\
&=&\sum_{j=0}^n\left( kB(t_{n+1},t_{j+\half})\eQ^{j+\half}
-\int_{t_j}^{t_{j+1}}B(t_{n+1},s)\eQ(s)\,ds,{\boldsymbol\chi} \right).
\end{eqnarray}
From the midpoint quadrature error, it follows that
$$
kg^{j+\half}-\int_{t_j}^{t_{j+1}}g(s)\,ds
= \int_{t_j}^{t_{j+\half}}(s-t_j)(s-t_{j+\half})D_s^2g(s)\,ds
+\int_{t_{j+\half}}^{t_{j+1}}(s-t_{j+1})(s-t_{j+\half})D_s^2g(s)\,ds.
$$
Use integration by parts to find that the boundary terms become zero and therefore, we arrive at 
\begin{eqnarray*}
kg^{j+\half}-\int_{t_j}^{t_{j+1}}g(s)\,ds&=&
-\int_{t_j}^{t_{j+\half}}\left(2s-(t_{j+\half}+t_j)\right) D_sg(s)\,ds
-\int_{t_{j+\half}}^{t_{j+1}}\left(2s-(t_{j+1}+t_{j+\half})\right) D_sg(s)\,ds
\end{eqnarray*}
and
$$
\left|kg^{j+\half}-\int_{t_j}^{t_{j+1}}g(s)\,ds\right|\leq 
\frac{k}{2}\int_{t_j}^{t_{j+1}}|D_sg(s)|\,ds,$$
Thus, (\ref{ss4}) with (\ref{nq}) leads to 
$$
\left| \ce^{n+1}(\eQ)({\boldsymbol\chi})\right| \leq 
Ck\sum_{j=0}^n\int_{t_j}^{t_{j+1}}\left(\|\eQ(s)\|+\|{\eQ}_t(s)\|\right)ds\,\|
{\boldsymbol\chi}\|,
$$
and 
\begin{eqnarray}\label{ss5} 
\|\ce^{n+1}(\eQ)\|&\leq& Ck\int_0^{t_{n+1}}\left(\|\eQ(s)\|+\|{\eQ}_t(s)\|\right)ds
\nonumber\\ &\leq& Ckht_{n+1}(\|u_0\|_3+\|u_1\|_2).
\end{eqnarray}
On substitution of (\ref{ss3}) and (\ref{ss5}) in (\ref{ss2}), it follows that
$$
k\sum_{n=0}^{m}\|\ce^{n+\half}(\tilde{\bq}_h)\|\leq C(T)(kh+k^2)
(\|u_0\|_3+\|u_1\|_2).
$$
Following similar line of proof, we can easily show that the last two terms in 
(\ref{ss1}) are also bounded by $C(T)(kh+k^2)(\|u_0\|_3+\|u_1\|_2)$.
Finally, by using the 
triangle inequality and the estimates in Lemmas~\ref{lm2} and \ref{lm3}
we complete the proof of the  theorem. $\hfill\Box$
%

\section {\bf Error Estimates for the Standard Mixed Method}
\se

Now, we extend our analysis to discuss optimal error estimates for 
$u$ and $\bs,$ satisfying the standard mixed method
(\ref{weak-mixed-1})-(\ref{weak-mixed-2}) with minimal regularity of the initial data. 
We first recall the following regularity results.
%
\begin{Lemma}\label{fap01}
 Let $(u,\bs)$ satisfy $(\ref{weak-mixed-1})$-$(\ref{weak-mixed-2})$. Then,
$$
\|D_t^ju(t)||+||D_t^{j-1}u(t)\|_1+\|D_t^{j-1}\bs(t)||\leq
C(T)(\|u_0\|_j+\|u_1\|_{j-1}),\quad j=1,\cdots,4,
$$
and
$$
\|D_t^ju(t)||_2\leq C(T)(\|u_0\|_{j+2}+\|u_1\|_{j+1}),\quad j=0,1,2.
$$ 
\end{Lemma}
With $W_h$ and $\bV_h$ defined as in Section 2, we define the corresponding
semidiscrete mixed finite element approximation to
(\ref{weak-mixed-1})-(\ref{weak-mixed-2}) as: Find a pair $(u_h,\bs_h)\in W_h\times \bV_h$
such that
\begin{eqnarray}
 (\alpha\bs_h,\bv_h)+\int_0^t (M(t,s)\bs_h(s),\bv_h)ds+(\nabla\cdot\bv_h,u_h)=0
 ~~~\forall ~\bv_h \in \bV_h \label{umfd}\\
 (u_{htt},w_h)-(\nabla\cdot\bs_h,w_h) = 0~~~~~~~\forall~w_h\in W_h,\label{bsmfd}
\end{eqnarray}
with $u_h(0)=P_h u_0$, and  $u_{ht}(0)=P_h u_1$.

Below, we present the main theorem of this section.
\begin{theorem}\label{mt03}
Let $(u,\bs)$ and $(u_h,\bs_h)$ satisfy  $(\ref{weak-mixed-1})$-$(\ref{weak-mixed-2})$ 
and $(\ref{umfd})$-$(\ref{bsmfd})$, respectively, with $u_h(0)=P_hu_0$ and 
$u_{ht}(0)=P_hu_1.$ Then, there exists a positive
constant $C,$ independent of $h,$ such that for $t\in (0,T]$
$$ \|u(t)-u_h(t)\| \le Ch^2\left(\|u_0\|_3+\|u_1\|_2\right). $$
\end{theorem}

\subsection {\bf Mixed Ritz-Volterra projections}

We discuss some of the properties of the mixed Ritz-Volterra 
projections used in our analysis. We formulate the mixed Ritz-Volterra projections as follows.
Given $(u(t),~\bs(t))\in W\times \bV,$ for $ t\in (0,T],$ find $(\tilde u_h,\tilde\bs_h):
(0,T]\longrightarrow W_h\times \bV_h $ satisfying 
\begin{eqnarray}
 (\alpha\es,\bv_h)+\int_0^t(M(t,s)\es(s),\bv_h)~ds+(\nabla\cdot\bv_h,\eu) =0,~~~
 \bv_h\in\bV_h, \label{eu01} \\
 (\nabla\cdot\es,w_h)=0,~~~w_h\in W_h, \label{eu02}
\end{eqnarray} 
where $\eu:= (u-\tilde u_h)$ and $\es:= (\bs-\tilde\bs_h).$


The following lemma can be easily obtained by combining Theorem 2.6 
of \cite{ELSWZ} and Lemma~\ref{fap01}.

\begin{Lemma}\label{lm62}
 Let $(\eu,\es)$ satisfy the system $(\ref{eu01})$-$(\ref{eu02})$. Then, there 
 is a positive constant C independent of $h$ such that
\begin{equation}\label{e1-n}
 \|D_t^j\es(t) \|\leq C h^r \left(\|u_0\|_{j+r+1}+\|u_1\|_{j+r}\right),\quad j=0,1,2,\quad  r=1,2,
\end{equation}
and
 \begin{equation}\label{e2-n}
 \|D_t^j\eu(t) \|\leq C h^2 \left(\|u_0\|_{j+2}+\|u_1\|_{j+1}\right),\quad j=0,1,2.
 \end{equation}
\end{Lemma}


\subsection {\bf Error Estimates}

Here, we discuss the proof of Theorem~\ref{mt03}. Set
$$ e_u:= u-u_h=\eu-\xu~~\mbox{ and }~~\be_{\bs}:= \bs-\bs_h=\es-\xs, $$
with $\xu=u_h-\tilde u_h$ and $\xs=\bs_h-\tilde\bs_h$.
Then, $(e_u,\be_{\bs})$ satisfy the following error equations
\begin{eqnarray} 
 (\alpha\be_{\bs},\bv_h)+\int_0^t (M(t,s)\be_{\bs}(s),\bv_h)ds+(\nabla\cdot\bv_h,e_u)=0
 ~~~\forall ~\bv_h \in \bV_h, \label{e5}\\
 (e_{u_{tt}},w_h)-(\nabla\cdot\be_{\bs},w_h) = 0~~~~~~~\forall~w_h\in W_h. \label{e4}
\end{eqnarray}
Using the mixed Ritz-Volterra projections, we write the above equations in terms of $\xu,\xs$ as
\begin{eqnarray}
 (\alpha\xs,\bv_h)+\int_0^t (M(t,s)\xs(s),\bv_h)ds+(\nabla\cdot\bv_h,\xu)=0,
 ~~~\forall ~\bv_h \in \bV_h, \label{x4}\\
 (\xu_{tt},w_h)-(\nabla\cdot\xs,w_h)= (\eu_{tt},w_h)~~~~~~~\forall~w_h\in W_h.\label{x5}
\end{eqnarray}
Below, we present a proof of our main theorem.

\noindent {\bf Proof of Theorem \ref{mt03}}: \\
Since the estimate $\eu$ is known from Lemma~\ref{lm62}, it is enough to estimate $\xu.$ 
We first integrate (\ref{x5}) and use the  fact that $u_{ht}(0)=P_hu_1$ to obtain
\begin{equation}\label{eq0}
(\xut,w_h)-(\nabla\cdot\bxs,w_h)= (\eu_{t},w_h)~~~~~~~\forall~w_h\in W_h.
\end{equation}
Now, choose $v_h=\bxs$ and $w_h=\xu$ in (\ref{x4}) and (\ref{eq0}), respectively, 
then add the resulting equations and use integration by parts to obtain
\begin{eqnarray*}
\frac{1}{2}\frac{d}{dt}\left[\|\xu\|^2+\|\alpha^\half\bxs\|^2\right]&=& 
(\eut,\xu)-\int_0^t (M(t,s)\xs(s),\bxs(t))\,ds\\
&=& (\eut,\xu)-(M(t,t)\bxs(t),\bxs(t))\,ds+\int_0^t (M_s(t,s)\bxs(s),\bxs(t))\,ds.
\end{eqnarray*}
Integrating from $0$ to $t$, and using the Cauchy-Schwarz inequality and 
the bounds for $M$, we immediately obtain
\begin{eqnarray*}
\|\xi_{u}(t)\|^2 &+& \|\alpha^\half\bxs(t)\|^2 \leq \|\xi_{u}(0)\|^2+2\int_0^t\|\eta_{u_{t}}(s)\|\|\xi_{u}(s)\|ds\\
&&+C(a_0,a_1,T)\left(\int_0^t\|\alpha^\half\bxs(s)\|^2\,ds+ 
\|\alpha^\half\bxs(t)\|\int_0^t\|\alpha^\half\bxs(s)\|\,ds\right). 
\end{eqnarray*}
Following the same arguments as in the proof of Theorem~\ref{thm2}, we deduce that 
$$
\|\xi_{u}(t)\|+\|\alpha^\half\bxs(t)\|
\leq C\left(\|\xi_{u}(0)\| +\int_0^t||\eta_{u_{t}}||\,ds\right).
$$
Finally, a use of the triangle inequality and Lemma~\ref{lm62} concludes the proof of 
theorem. \hfill{$\Box$}

The discrete-in-time scheme for  problem
$(\ref{weak-mixed-1})$-$(\ref{weak-mixed-2})$ is to seek  
$(U^{n},\Z^{n})\in W_h\times{\bV}_h$, such that
\begin{eqnarray}
&&\frac{2}{\dt}(\ph U^{\half},w_h)-(\nabla\cdot\Z^{\half},w_h)=(\frac{2}{\dt}u_1,w_h), \label{7-1-s}\\
&&(\alpha\Z^{n+\half},{\bv}_h)+(U^{n+\half},\nabla\cdot\bv_h)+
\ii^{n+\half}({\cal M}^{n+\half}(\Z,{\bv}_h)=0,\; n\geq 0,\label{7-1b-s}\\
&&(\ps U^{n},w_h)-(\nabla\cdot\Z^{n;\qr},w_h)=0,\; n\geq 1,\label{7-1c-s}
\end{eqnarray}
 for all $(w_h,{\bv}_h)\in W_h\times{\bV}_h,$ with
given $(U^{0},\Z^{0})\in W_h\times\bV_h$. In (\ref{7-1b-s}),
$$
\ii ^{n+\half}({\cal M}^{n+\half}(\Z,{\bv}_h))=\frac{1}{2}\left(\ii^{n+1}
({\cal M}^{n+1}(\Z,{\bv}_h)+\ii^n(\B^n(\Z,{\bv}_h)) \right),
$$
where
$$\ii^n({\cal M}^n(\Z,\chi))=\dt\sum_{j=0}^{n-1}(M(t_n,t_{j+1/2})\Z^{j+1/2},\chi).$$
This choice of the time  discretization leads to a second order accuracy in $\dt$.
Below we state the following theorem. Since the proof follows exactly the steps leading to
Theorem~\ref{thm4} with appropriate changes, we skip the proof.
\begin{theorem}\label{thm-7-2}
Let $(u,\bs)$ be the solution of $(\ref{weak-mixed-1})$-$(\ref{weak-mixed-2})$
and $(U^n,\Z^n)\in W_h\times{\bV}_h$ be the 
solution of $(\ref{7-1-s})$-$(\ref{7-1c-s})$.
Assume that $U^0=P_hu_0$ and $\Z^0=P_h(A\nabla u_0)$.
Then, there exists a positive constant $C$, independent of $h$ and $k$,
such that for small $k$ with $k=O(h)$, 
\begin{equation}\label{6-2-bb-n}
\|u(t_{m+1})-U^{m+1}\|
\leq  C (h^2+k^2)\left(\|u_0\|_3+\|u_1\|_2\right), \quad m=0,1,\cdots,N-1.
\end{equation}
\end{theorem}

{\bf { Acknowledgements.}}  
The two authors gratefully acknowledge the research support of the Department of Science and Technology, 
Government of India through the National Programme on Differential Equations: Theory, Computation and 
Applications vide DST Project No.SERB/F/1279/2011-2012, and 
the support by Sultan Qaboos University under Grant IG/SCI/DOMS/13/02.

\end{document}